\documentclass[12pt]{article}

\usepackage[english]{babel}
\usepackage{graphicx}
\usepackage{framed,caption}
\usepackage[normalem]{ulem}
\usepackage{indentfirst}
\usepackage{amsmath,amsthm,amssymb,amsfonts}
\usepackage[italicdiff]{physics}
\usepackage[T1]{fontenc}

\usepackage{setspace}
\usepackage{wrapfig}
\usepackage{lmodern,mathrsfs}
\usepackage[inline,shortlabels]{enumitem}
\setlist{topsep=2pt,itemsep=2pt,parsep=3pt,partopsep=2pt}

\usepackage[dvipsnames]{xcolor}
\usepackage[utf8]{inputenc}
\usepackage[a4paper, top=0.8in,bottom=0.8in, left=0.8in, right=0.8in, footskip=0.3in, includefoot]{geometry}
\usepackage[most]{tcolorbox}
\usepackage{tikz,tikz-3dplot,tikz-cd,tkz-tab,tkz-euclide,pgf,pgfplots}
\pgfplotsset{compat=newest}
\usepackage{multicol}
\usepackage[bottom,multiple]{footmisc} 

\usepackage[backend=bibtex,style=numeric]{biblatex}
\usepackage{hyperref}
\usepackage[nameinlink]{cleveref} 

\usepackage{graphicx}
 \graphicspath{{figures/}} 

\numberwithin{equation}{section}

\usepackage{framed,color}
\definecolor{shadecolor}{rgb}{0.9412,0.9725, 1}
\usepackage{caption}
\usepackage{tcolorbox}

\usepackage{multirow}
\usepackage{epstopdf}
\usepackage{listings}
\usepackage[]{mcode}
\usepackage{algorithm}
\usepackage{latexsym}
\usepackage{showidx}
\usepackage{latexsym}
\usepackage{amssymb}

\newcommand{\ignore}[1]{}
\usepackage{todonotes}

\usepackage{tikz}

\newcommand{\hide}[1]{} 

\renewcommand{\qed}{\hfill\blacksquare}

\newcommand{\ep}{\varepsilon}

\newcommand{\C}{\mathbb{C}}

\newcommand{\F}{\mathbb{F}}

\newcommand{\N}{\mathbb{N}}

\newcommand{\R}{\mathbb{R}}

\newcommand{\fl}{\mathrm{fl}}

\newcommand{\wh}{\widehat}
\newcommand{\xh}{\widehat{x}}
\newcommand{\yh}{\widehat{y}}

\newcommand{\sh}{\widehat{s}}

\renewcommand{\iff}{\Leftrightarrow}

\newtheoremstyle{mystyle}{}{}{}{}{\sffamily\bfseries}{.}{ }{}
\newtheoremstyle{cstyle}{}{}{}{}{\sffamily\bfseries}{.}{ }{\thmnote{#3}}
\makeatletter
\renewenvironment{proof}[1][\proofname] {\par\pushQED{\qed}{\normalfont\sffamily\bfseries\topsep6\p@\@plus6\p@\relax #1\@addpunct{.} }}{\popQED\endtrivlist\@endpefalse}
\makeatother
\theoremstyle{mystyle}{\newtheorem{definition}{Definition}[section]}
\theoremstyle{mystyle}{\newtheorem{theorem}[definition]{Theorem}}
\theoremstyle{mystyle}{\newtheorem{lemma}[definition]{Lemma}}
\theoremstyle{mystyle}{}
\theoremstyle{mystyle}{}
\theoremstyle{mystyle}{}
\theoremstyle{mystyle}{}
\theoremstyle{mystyle}{\newtheorem{workout}[definition]{Workout}}
\theoremstyle{mystyle}{}

\usepackage{amsmath}
\newtheorem{remark}{Remark}[section]
\newtheorem{example}{Example}[section]

\theoremstyle{cstyle}{}

\newtheoremstyle{warn}{}{}{}{}{\normalfont}{}{ }{}
\theoremstyle{warn}

\newcommand{\warningsign}[1]{\tikz[scale=#1,every node/.style={transform shape}]{\draw[-,line width={#1*0.8mm},red,fill=yellow,rounded corners={#1*2.5mm}] (0,0)--(1,{-sqrt(3)})--(-1,{-sqrt(3)})--cycle;
\node at (0,-1) {\fontsize{48}{60}\selectfont\bfseries!};}}

\tcolorboxenvironment{definition}{boxrule=0pt,boxsep=0pt,colback={red!10},left=8pt,right=8pt,enhanced jigsaw, borderline west={2pt}{0pt}{red},sharp corners,before skip=10pt,after skip=10pt,breakable}
\tcolorboxenvironment{proposition}{boxrule=0pt,boxsep=0pt,colback={Orange!10},left=8pt,right=8pt,enhanced jigsaw, borderline west={2pt}{0pt}{Orange},sharp corners,before skip=10pt,after skip=10pt,breakable}
\tcolorboxenvironment{theorem}{boxrule=0pt,boxsep=0pt,colback={blue!10},left=8pt,right=8pt,enhanced jigsaw, borderline west={2pt}{0pt}{blue},sharp corners,before skip=10pt,after skip=10pt,breakable}
\tcolorboxenvironment{lemma}{boxrule=0pt,boxsep=0pt,colback={blue!10},left=8pt,right=8pt,enhanced jigsaw, borderline west={2pt}{0pt}{blue},sharp corners,before skip=10pt,after skip=10pt,breakable}
\tcolorboxenvironment{corollary}{boxrule=0pt,boxsep=0pt,colback={violet!10},left=8pt,right=8pt,enhanced jigsaw, borderline west={2pt}{0pt}{violet},sharp corners,before skip=10pt,after skip=10pt,breakable}
\tcolorboxenvironment{proof}{boxrule=0pt,boxsep=0pt,blanker,borderline west={2pt}{0pt}{CadetBlue!80!white},left=8pt,right=8pt,sharp corners,before skip=10pt,after skip=10pt,breakable}
\tcolorboxenvironment{remark}{boxrule=0pt,boxsep=0pt,blanker,borderline west={2pt}{0pt}{Green},left=8pt,right=8pt,before skip=10pt,after skip=10pt,breakable}
\tcolorboxenvironment{remarks}{boxrule=0pt,boxsep=0pt,blanker,borderline west={2pt}{0pt}{Green},left=8pt,right=8pt,before skip=10pt,after skip=10pt,breakable}
\tcolorboxenvironment{example}{boxrule=0pt,boxsep=0pt,blanker,borderline west={2pt}{0pt}{Black},left=8pt,right=8pt,sharp corners,before skip=10pt,after skip=10pt,breakable}
\tcolorboxenvironment{examples}{boxrule=0pt,boxsep=0pt,blanker,borderline west={2pt}{0pt}{Black},left=8pt,right=8pt,sharp corners,before skip=10pt,after skip=10pt,breakable}
\tcolorboxenvironment{cthm}{boxrule=0pt,boxsep=0pt,colback={gray!10},left=8pt,right=8pt,enhanced jigsaw, borderline west={2pt}{0pt}{gray},sharp corners,before skip=10pt,after skip=10pt,breakable}
\tcolorboxenvironment{workout}{boxrule=0pt,boxsep=0pt,colback={Cyan!10},left=8pt,right=8pt,enhanced jigsaw, borderline west={2pt}{0pt}{Cyan},sharp corners,before skip=10pt,after skip=10pt,breakable}
\tcolorboxenvironment{labexercise}{boxrule=0pt,boxsep=0pt,colback={green!10},left=8pt,right=8pt,enhanced jigsaw, borderline west={2pt}{0pt}{green},sharp corners,before skip=10pt,after skip=10pt,breakable}

\newenvironment{talign*}{\let\displaystyle\textstyle\csname align*\endcsname}{\endalign}

\usepackage[explicit]{titlesec}
\titleformat{\section}{\fontsize{18}{30}\sffamily\bfseries}{\thesection}{18pt}{#1}
\titleformat{\subsection}{\fontsize{15}{18}\sffamily\bfseries}{\thesubsection}{15pt}{#1}
\titleformat{\subsubsection}{\fontsize{10}{12}\sffamily\large\bfseries}{\thesubsubsection}{8pt}{#1}

\titlespacing*{\section}{0pt}{5pt}{5pt}
\titlespacing*{\subsection}{0pt}{5pt}{5pt}
\titlespacing*{\subsubsection}{0pt}{5pt}{5pt}

\DeclareMathAlphabet\mathbfcal{OMS}{cmsy}{b}{n}
\newcommand{\cond}{\mathrm{cond}}
\newcommand{\ds}{\displaystyle}
\newcommand{\fs}{\footnotesize}
\newcommand{\vsp}{\vspace{0.3cm}}

\title{\large\sffamily\bfseries Lecture Notes \\\LARGE{Fundamentals of Computing}}
\author{\Large\sffamily Davoud Mirzaei\\ Uppsala University}
\date{\sffamily June 10, 2024}

\begin{document}

\setlength{\abovedisplayskip}{0pt}
\setlength{\belowdisplayskip}{3pt}
\setlength{\abovedisplayshortskip}{0pt}
\setlength{\belowdisplayshortskip}{0pt}
\maketitle

\definecolor{contcol1}{HTML}{72E094}
\definecolor{contcol2}{HTML}{24E2D6}
\definecolor{convcol1}{HTML}{C0392B}
\definecolor{convcol2}{HTML}{8E44AD}

\begin{tcolorbox}[title=Contents, fonttitle=\huge\sffamily\bfseries\selectfont,interior style={left color=contcol1!10!white,right color=contcol2!10!white},frame style={left color=contcol1!30!white,right color=contcol2!30!white},coltitle=black,top=2mm,bottom=2mm,left=2mm,right=2mm,drop fuzzy shadow,enhanced,breakable]
\makeatletter
\@starttoc{toc}
\makeatother
\end{tcolorbox}

\vspace*{10mm}

\thispagestyle{empty}
\newpage
\pagenumbering{arabic}

This lecture addresses some general ideas behind numerical
computations ranging from representation of numbers in computers to stability and accuracy of standard algorithms for some simple mathematical problems. Some parts of the lecture follow \cite{Dahlquist-Bjork:2008}, \cite{Higham:2002}, and \cite{Gautschi:2012}.


\section{What is scientific computing?}
\begin{figure}[!th]
\centering
\includegraphics[scale=0.35]{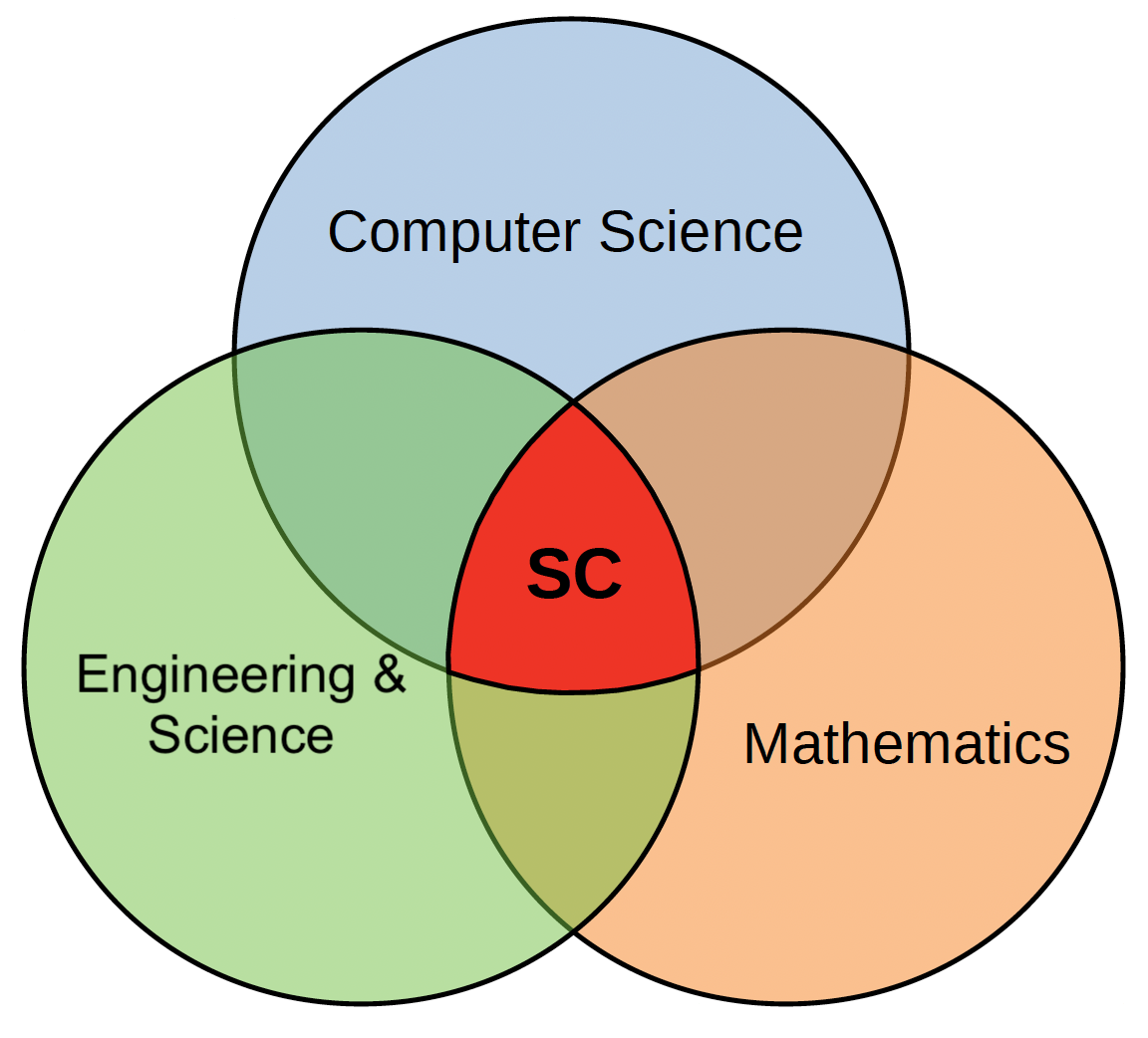}
\caption{The land of scientific computing (SC)}
\label{fig:sci_comput}
\end{figure}

{\em Scientific computing} is concerned with the design and analysis of algorithms for solving mathematical
problems that arise in science and engineering.
It is distinguished from most other parts of computer science in that it deals with quantities that are continuous, as opposed to discrete. Continuous quantities are functions and equations whose underlying variables (time, distance,
velocity, temperature, density, pressure, stress, ...) are continuous in nature.
This subject is also called  {\em numerical analysis} or {\em computational mathematics} but nowadays it is mostly refereed to as
{\em scientific computing}. Figure \ref{fig:sci_comput} shows that scientific computing could be viewed as the intersection of
computer science, applied mathematics and science and engineering.

Most problems in continuous mathematics, such as those involving derivatives, integrals, or nonlinearities, cannot be solved exactly (analytically) and must be addressed using approximate processes that eventually converge to numerical solutions. A key aspect of scientific computing is the development of convergent algorithms and the analysis of the accuracy of these approximations. Therefore, another critical factor in scientific computing is its focus on the effect of {\em approximations}. We summarize that the
key distinguishing features of scientific computing include
\begin{itemize}
\item dealing with continuous quantities (e.g., time, distance, velocity, temperature, density, pressure) typically measured by real numbers, and
\item considering the effects of approximations.
\end{itemize}

\subsection{Computational simulation}

The classic pair of opposed but mutually supporting scientific paradigms are {\em theory} and {\em experimentation}. A third paradigm, {\em computational simulation}, emerged through the work of John von Neumann and others in the mid-20th century\footnote{Jim Gray, a 1998 Turing Award winner and a leading computer scientist, proposed a ``fourth paradigm'' in scientific research in one of his last talks in 2007. This paradigm, data-intensive science, is a methodological approach to discovery based on data analysis, extending beyond theoretical and experimental research and computational simulations.}. Computational simulation involves representing and emulating physical systems or processes using computers.

Today, computation has become an equal and indispensable partner alongside theory and experiment in the pursuit of knowledge and technological advancement. Numerical simulation allows the study of complex systems and natural phenomena that would be too expensive, dangerous, or even impossible to study theoretically or through direct experimentation.
For instance, in astrophysics, the detailed behavior of two colliding black holes is too complex to determine theoretically and impossible to observe directly or replicate in a laboratory. To simulate this scenario computationally, one needs an appropriate mathematical representation (such as Einstein's equations of general relativity), an algorithm to solve these equations numerically, and a sufficiently powerful computer to implement the algorithm. Similarly, in improving automobile safety, crash testing on a computer is far less expensive and dangerous than real-life testing. This allows for a more thorough exploration of potential design parameters, leading to the development of optimal designs.

The overall problem-solving process in computational simulation typically includes the following steps:
\begin{enumerate}
\item 
{\bf Develop a mathematical model:} This involves formulating a mathematical representation, usually in the form of equations, of the physical phenomenon or system of interest.
\item 
{\bf Develop algorithms:} Create numerical algorithms to solve these equations.
\item 
{\bf Implement and execute the algorithms:} Convert the algorithms into computer software and run the simulations.
\item 
{\bf Interpret and validate the results:} Analyze the computed results and verify their accuracy, repeating any or all of the preceding steps as necessary.
\end{enumerate}
Step 1, known as mathematical modeling, requires specialized knowledge of the relevant scientific or engineering disciplines, as well as applied mathematics. Steps 2 and 3, which involve designing, analyzing, implementing, and using numerical algorithms and software, constitute the core focus of scientific computing.

\begin{figure}
\begin{center}
\begin{tabular}{ccc}
\begin{tabular}{l}
\\ 
\includegraphics[scale=.065]{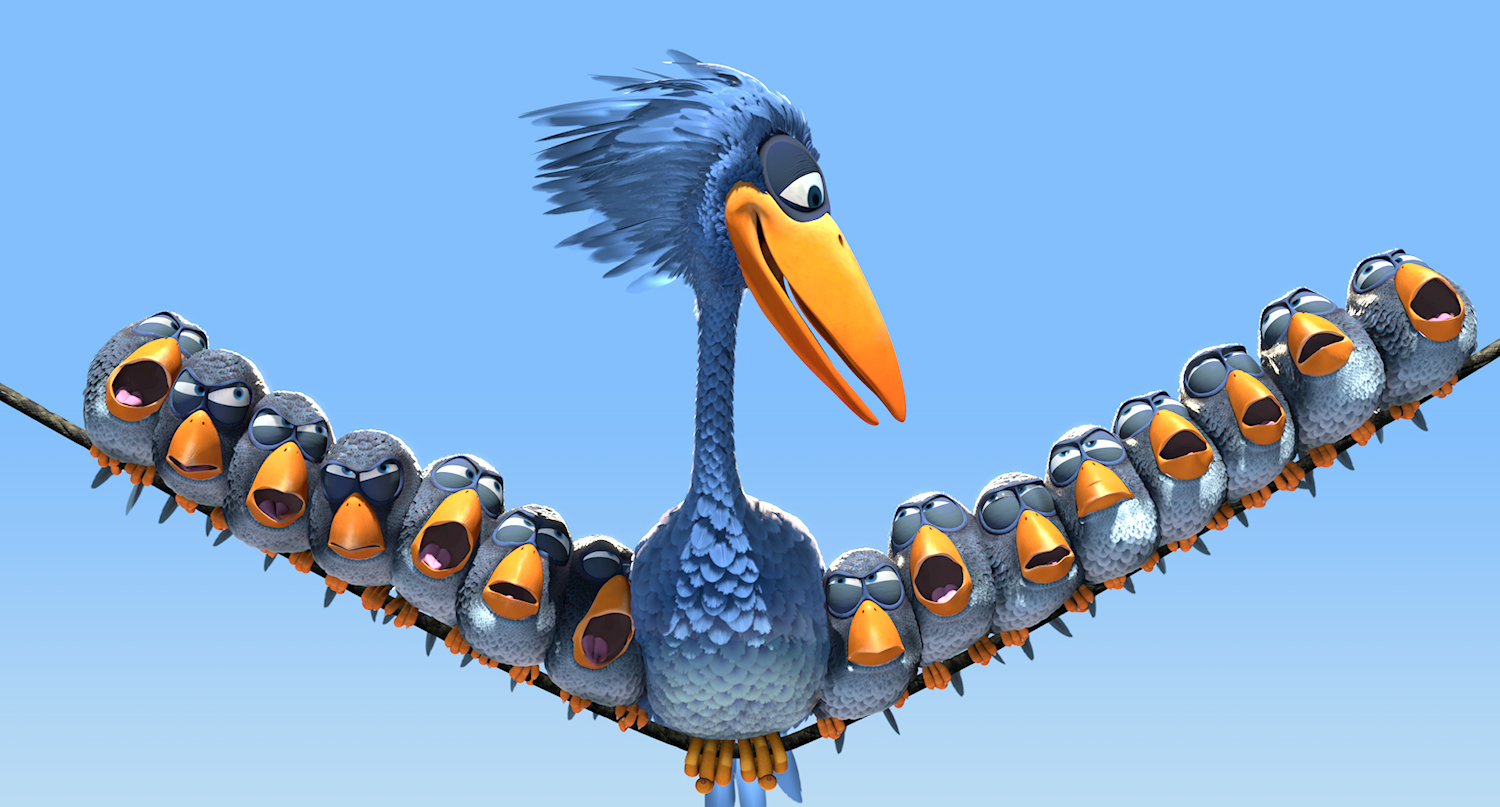}
\end{tabular}
 & 
\tikz [line width=.7mm]
\draw [color = red,arrows = {-Computer Modern Rightarrow[round]}] (0,0) -- (2.8,0) node[midway,above] {modelling}; 
 &
 \hspace{-.7cm}
\colorbox{green!50}{
\begin{tabular}{l}
$-y''(x) = f(x), \; 0\leqslant x\leqslant \ell $\\
$y(0)=y(\ell) = 0$ 
\end{tabular}}\\
\tikz [line width=.7mm]
\draw [color = red,arrows = {-Computer Modern Rightarrow[round]}] (0,0) -- (0,1.2) node[midway,left] {interpretation}; \qquad 
 &   &   
\tikz [line width=.7mm]
\draw [color = red,arrows = {-Computer Modern Rightarrow[round]}] (0,0) -- (0,-1.2) node[midway,right] {discretization};  \\
\begin{tabular}{l}
\\ 
\includegraphics[scale=.20]{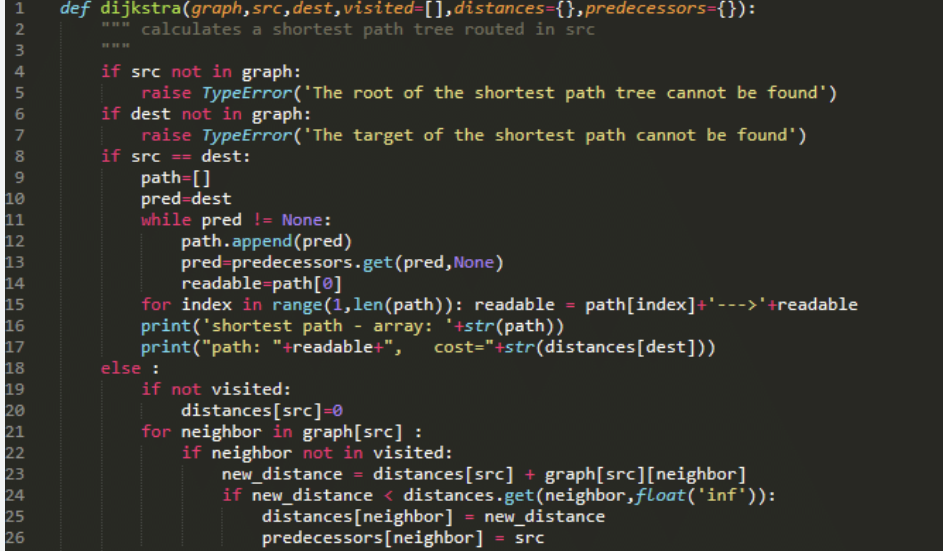}
\end{tabular}
  &   
\hspace{-.3cm}
\tikz [line width=.7mm]
\draw [color = red,arrows = {-Computer Modern Rightarrow[round]}] (0,0) -- (-2.8,0) node[midway,above] {code implement.}; 
  &   
  \hspace{-.5cm}
\colorbox{yellow!50}{
\begin{tabular}{l}
$y_{k+1}-2y_k+y_{k-1}=-h^2f_k $\\
$y_0=y_N = 0$ 
\end{tabular}}
\end{tabular}
\caption{Problem solving steps: birds on the wire!}\label{fig_birds}
\end{center}
\end{figure}

As a simple and funny example, consider a bird family (mother bird and chicks) sitting on an elastic wire of length \(\ell\) meters, fixed at both ends. Under some mild simplifying assumptions, the displacement of the wire from the horizontal line (denoted as \(y\)) satisfies a simple second-order ordinary differential equation (ODE) with boundary conditions (see Figure \ref{fig_birds}). 
To solve this, we can use the finite difference method (FDM) to discretize the underlying ODE and boundary conditions. In the third step, we implement the discretized equations in a computer program to obtain a numerical (approximate) solution, \(y = (y_0, y_1, \ldots, y_N)\), at discrete points \(x_0, x_1, \ldots, x_N\) in interval \([0, \ell]\). This numerical solution approximates the exact continuous solution \(y(x)\) for \(x \in [0, \ell]\).
Finally, we interpret and validate our numerical solution by comparing it to the real phenomenon, i.e., the birds on the wire.

As a real-world problem, consider the simulation of the aerodynamics of a car to enhance its efficiency and reduce fuel consumption. This problem involves the mathematical modeling of two fundamental physical phenomena, fluid dynamics and solid mechanics. Specifically, it requires solving the Navier-Stokes equations for fluid flow (air) and the elasticity equations for the solid structure (car body). These equations are partial differential equations (PDEs) that describe the conservation of physical quantities such as momentum and energy.
In practical scenarios, finding exact solutions to these PDEs using analytical techniques is often impossible due to their complexity. This is where computational methods come into play. For the fluid domain, the Finite Volume Method (FVM) is commonly employed, while the Finite Element Method (FEM) is typically used for the solid domain. To ensure the reliability of these numerical solutions, it is essential to analyze the algorithms in terms of stability, convergence, and computational complexity. Stability analysis ensures that the numerical solution behaves correctly as it progresses, while convergence analysis guarantees that the solution approximates the true solution as the computational parameters are refined. Computational complexity assesses the efficiency of the algorithm in terms of time and resources required.

\begin{figure}[ht!]
\begin{center}
\includegraphics[scale=.3]{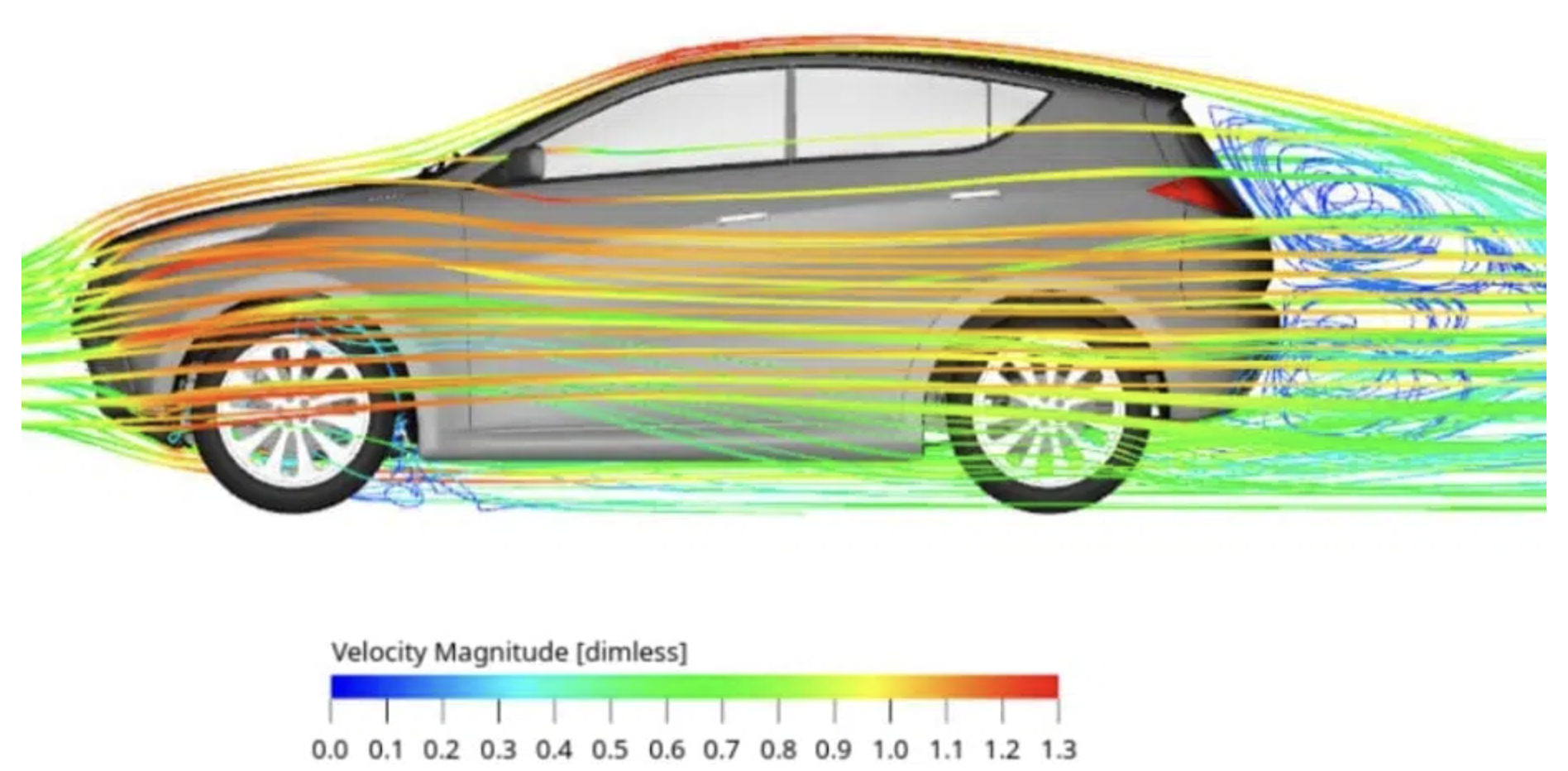}
\caption{Simulation of the velocity filed and pressure of the air surrounding a car (image from \texttt{www.vias3d.com})}\label{fig_car_sim}
\end{center}
\end{figure}

Addressing such problems requires a deep understanding of mathematical theory, proficiency in computational techniques, and expertise in high-performance computing. These skills enable us to develop, implement, and optimize algorithms that deliver accurate and efficient solutions to complex scientific and engineering problems.
Figure \ref{fig_car_sim} illustrates an example of such a simulation, showing the velocity field and its magnitude (represented by colors) for a sample car. This example not only demonstrates the practical application of numerical methods but also highlights the importance of interdisciplinary knowledge in scientific computing.  
\vsp

\section{Sources of approximation}

In the course of a numerical algorithm, computational results are influenced by various types of approximation errors. These errors can propagate from their sources to subsequently computed quantities, sometimes with significant amplification or damping.
Some approximations may occur before a computation begins:
\begin{itemize}
\item[a.] \textbf{Simplifications in the mathematical model.}
When developing a mathematical model for a natural phenomenon, certain simplifying assumptions are often made to facilitate the modeling process. These assumptions help reduce complexity but may introduce minor deviations from real-world behavior. For example,
in modeling a pendulum, we might assume the string is massless and ignore air resistance. In the heat conduction problem of a rod, we often assume the rod is composed of a homogeneous material. In economic calculations, we might assume that the interest rate remains constant over a given period. 
These assumptions are essential for making complex models manageable and solvable, but they do so at the cost of introducing small deviations from the actual physical solutions.

\item[b.] \textbf{Errors in input data.} Input data can often be the result of measurements that have been contaminated by various types of errors. The finite precision of laboratory instruments introduces systematic errors, while variations in the experimental environment lead to random errors in the input data. Additionally, input data may have been generated by a previous computational step, the results of which were only approximate.
\end{itemize}
These types of approximation errors are generally considered uncontrollable in numerical analysis. However, providing feedback to the developer of the mathematical model can sometimes be beneficial. Such communication can lead to refinements in the model and reduce the impact of these errors on the final results.
In scientific computing, however, the focus is mostly on two other types of errors:
\begin{itemize}
\item[c.] \textbf{Discretization errors} (sometimes calles \textbf{ truncation error}),
\item[d.] \textbf{Rounding errors} (also called \textbf{roundoff errors}).
\end{itemize}
The main distinction between rounding and discretization errors is that rounding errors arise from {\em arithmetic} calculations, i.e., manipulation of {\em digits} of numbers (so they mainly have a computer arithmetic nature) while discretization errors stem from {\em algorithmic} calculations (so they mainly have a {\em mathematical} nature). Although the most parts of this lecture are devoted to understanding machine arithmetic and rounding errors, here we categorize some sources of discretization errors:
\begin{itemize}
\item[(a)] replacing an infinite process by a finite approximation, for example,
\begin{itemize}
\item (discrete case) replacing an infinite series by a summation of a finite number of terms
\item (continuous case) replacing an integral of a function by a finite summation of values of the function as in the trapezoidal or Simpson's rules.
\end{itemize}

\item[(b)] replacing an infinitesimal process by a finite approximation, e.g., replacing the limit in differentiation by a finite approximation

\item[(c)]truncation of an iteration that, in theory, should continue forever after a finite number of iterations in practice. This happens e.g., in iterative methods for solving
\begin{itemize}
\item nonlinear equations (rootfinding) like the Newton-Raphson method
\item linear {\em systems} of equations like the Jacobi or Gauss-Seidel, SOR, conjugate gradient, etc.
\end{itemize}
Doing these, the contribution of the remaining terms or iterations are not taken into account. This is why this type of error is also called the `truncation error'.
\end{itemize}

\vsp
\begin{example}
The surface area of the Earth cab be computed
using the formula
$$
A = 4\pi r^2
$$
for the surface area of a sphere of radius $r$. The use of this formula for the computation
involves a number of approximations:
\begin{itemize}
\item The Earth is modeled as a sphere, which is an idealization of its true shape.
\begin{center}
\includegraphics[scale=.25]{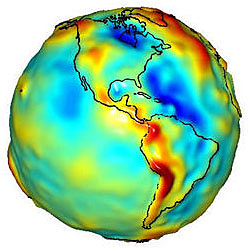} \includegraphics[scale=.6]{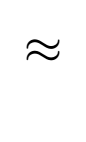}   \includegraphics[scale=.24]{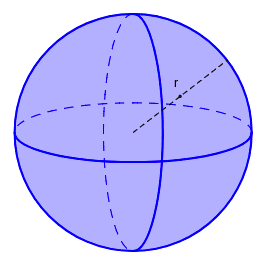}
\end{center}
\item The value for the radius, $r \approx 6370$ km, is based on a combination of empirical
measurements and previous computations.
\item The value for $\pi$ is given by an infinite limiting process, which must be truncated
at some point, e.g. $\pi\approx 3.1415$.
\item The numerical values for the input data, as well as the results of the arithmetic
operations performed on them, are rounded in a computer or calculator.
\end{itemize}
The accuracy of the computed result $A \approx 2.5146\times 10^{5} \mathrm{km}^2$ depends on all of these approximations.
\end{example}
\vsp
In the above example we observed the effects of modeling, input data and rounding errors.
To give an example for discretization error, let us recall the Taylor's expansion theorem which is one of the most fundamental results throughout mathematics whose importance is well beyond simple understanding of the discretization errors\footnote{Even in this course we will use Taylor's theorem later e.g., in the analysis of the discretization errors in derivation of numerical differentiation formulas and in the numerical solution of ODEs.}.

\begin{theorem}
Let $f$ be continuously differentiable up to order $n+1$ on the interval $[a, b]$ and $x_0 \in [a, b]$. Then, for every $x \in [a, b]$, there exists a point $\xi$ between $x_0$ and $x$ such that
\[
f(x) = p_n(x) + r_n(x),
\]
where
\begin{align*}
p_n(x) &= f(x_0) + f'(x_0) (x-x_0) + \frac{f''(x_0)}{2!} (x-x_0)^2 + \dots + \frac{f^{(n)}(x_0)}{n!} (x-x_0)^n,
\end{align*}
and
\[
r_n(x) = \frac{f^{(n+1)}(\xi)}{(n+1)!} (x-x_0)^{n+1}.
\]
\end{theorem}

Here, \(p_n\) is referred to as the Taylor polynomial of degree \(n\) of function \(f\) around \(x_0\), and \(r_n\) is the remainder associated with the polynomial \(p_n\). Since \(p_n(x)\) retains only the first \(n+1\) terms of the infinite series, it is appropriate to view the remainder \(r_n(x)\) as the discretization or truncation error corresponding to \(p_n(x)\). 
\vsp

\begin{example}
Assume that $f(x) = \exp(x)$ and $x_0=0$. Since all derivatives of exponential function at zero are $1$, the Taylor series expansion reads as
\vsp
\begin{equation}
\label{expSer:eq}
\exp(x) = 
\underbrace{1 + x + \frac{x^2}{2!} \dots + \frac{x^n}{n!}}_{p_n(x)} +
\underbrace{\frac{\exp(\xi)}{(n+1)!}x^{n+1}}_{r_n(x)},
\end{equation}
for a value $\xi$ between $x$ and $0$. 

Now, we use $n=3$ and compute an approximation for
$\sqrt{e}$ with an upper bound on the (truncation) error. 
We set $x = \frac{1}{2}$ and $n=3$ in \eqref{expSer:eq} to get
\[
p_3(1/2) = 1 + \frac{1}{2} + \frac{1}{8} +  \frac{1}{48} = \frac{79}{48}
\]
as an approximation for $\sqrt{e}$. The corresponding discretization error is
\[
r_3(1/2) = \frac{\exp(\xi)}{4!} (\frac{1}{2})^4 =  \frac{\exp(\xi)}{384},
\]
where $\xi \in (0, \frac{1}{2})$. Since exponential is an increasing function we have $\exp(\xi) < \exp(1/2) < 2$. Therefore, {\em an} upper bound for absolute value of discretization error is
 \[
  |r_3(1/2)| < \frac{2}{384} \approx 0.0052.
 \]
\end{example}
Let us review what we did in the above example. The problem was to find an approximation for $\sqrt{e}$. The algorithm we used to solve the problem was the degree-$3$ Taylor polynomial of the exponential function. The discretization error was the result of replacing an infinite number of terms in that series with a polynomial of degree three only.
In addition, note that we have {\em not} committed any rounding errors as long as we keep our approximation in the rational form $\frac{79}{48}$; the only error till now is the discretization error. However, rounding errors occur the moment we represent $\frac{79}{48}$ approximately e.g., as $1.6458$.
It can be shown that the first few digits of $\sqrt{e}$ are $1.648721270700128$. So, we observe that the first three digits $1.64$ of our approximation were indeed correct and that the amount of discretization error was indeed about
 $0.0029$ which is, of course, smaller than the upper bound we computed.
\vsp

\subsection{Absolute and relative errors}
Suppose that $\wh x$ is an approximation for a real number $x$. Two ways to measure the amount of error in $\wh x$ are the absolute error
\[
| x - \wh x|
\]
and the relative error defined as
\begin{equation}
\label{relErrorDef:eq}
\frac{| x - \wh x|}{|x|}
\vsp
\end{equation}
provided that $x \neq 0$. If we have both small and large quantities in a computation, it is the relative error that is more useful.

Since in most times the correct value $x$ is not known, we cannot compute the exact amount of (absolute or relative) error. Instead, we either approximate the error or try to find a bound for it. For instance, in Example 1.5 we obtained an upper bound for the amount of (absolute) discretization error in approximating $\sqrt{e}$ by $\frac{79}{48}$. In a similar spirit, we will see later in the fundamental theorem of rounding errors that we can find an upper bound for the amount of (relative) error resulting from rounding a real number to a machine number.

If $x$ is a vector then the absolute and the relative error are defined as
\vsp
$$
\|x-\wh{x}\| \quad \mbox{and} \quad \frac{\|x-\wh{x}\|}{\|x\|}
\vsp
$$
respectively, where $\|\cdot\|$ denotes some vector norm. For functions, norms defined on function spaces can be replaced.

 \section{Computer representation of numbers}
To analyse the algorithms, it is essential to understand how machines store numbers, represent them, and perform computations. In everyday life, we use base-10 numbers, likely because humans historically began counting with their fingers. Computers, however, use base-2 representation and binary arithmetic, as electrical devices typically distinguish between two states, e.g. lamps being off or on, or magnetic fields being clockwise or counterclockwise. However, here we discuss the representation of numbers in an arbitrary integer base $\beta \geq 2$.
\vsp 
\subsection{Fixed-point representation of numbers}
Given an integer base $\beta \geq 2$, for any real number $x$ there exists $n\in \N$ such that $x$ can be represented as a possibly infinite string like
\vsp
\begin{equation*}
\begin{split}
x & =  (-1)^{\sigma} (d_n \beta^n +  d_{n-1} \beta^{n-1} + \cdots + d_0 + d_{-1} \beta^{-1} + d_{-2} \beta^{-2} + \cdots) \\
  & =: \pm (d_n d_{n-1} \cdots d_0 . d_{-1} d_{-2} \cdots)_{\beta}
  \vsp
\end{split}
\end{equation*}
where $\sigma \in \{ 0, 1\}$ characterizes the sign of $x$ and the {\em digits} $d_k$ are integers in $\{0,1,\ldots, \beta-1\}$.
This is the \textbf{position system} in which one can give simple and general rules for the arithmetic operations. 
In computers, real or integer numbers are typically stored using $32$ or $64$ bits, which are referred to as {\em word lengths}. In the first generation of computers, calculations were performed using a {\bf fixed-point} number system. In this system, numbers have a fixed number of $p$ digits in the fractional part, i.e.,
\begin{equation*}\label{realNum:eq}
\wh x = \pm (d_n d_{n-1} \cdots d_0 . d_{-1} d_{-2} \cdots d_{-p})_{\beta}.
\end{equation*}
Let us denote the set of all representable numbers in this system by $F(\beta, n, p)$. 
The smaller the base is, the simpler these rules become. This is just one reason why most computers operate in base $2$, the {\em binary} number
system, where $d_k\in\{0,1\}$. The addition and multiplication then take the following simple form:
$$
0+0 = 0,\quad 0+1=1+0=1,\quad 1+1 = 10,\quad 0\cdot 0=0,\quad 0\cdot 1 = 1\cdot 0 =0,\quad 1\cdot 1 = 1.
$$
The digits in the binary system are called \textbf{bits} (\textbf{b}inary dig\textbf{its}).
If the computer's word length is \( s+1 \) bits (including the sign bit), then for \(\beta = 2\), the range of representable numbers is \( F(2, n, p) \subset [-2^{s-p}, 2^{s-p}] \). For example, if \( s+1 = 64 \) and \( n = p = 31 \), then \( 2^{s-p} \approx 4.2950 \times 10^9 \). This interval is not large enough for many real-world applications. Scientists work with numbers of vastly different scales, from biologists dealing with microscopic quantities to astrophysicists measuring immense distances between celestial bodies. Thus, it is crucial to have a number system capable of handling both very small and very large numbers.

The fixed-point number system provides the same level of sensitivity to both small and large numbers. This sacrifices the efficiency of this system. For instance, the distance between 0 and its next representable number in $F(\beta,n,p)$ is identical with the distance between the largest and second-largest numbers in this system. This wastes many bits for representing unnecessary numbers that could be instead used to accommodate a wider range of numbers. 

\begin{workout}
What is the number of numbers in $F(\beta,n,p)$. Determine the largest and smallest positive numbers in this system. Show the distance between two consecutive numbers in $F(\beta,n,p)$ is $\beta^{-p}$ (the distance is independent of the magnitude of numbers).
\end{workout}
\vsp 

\subsection{Floating-point representation of numbers}

Every real number $x$ can be represented as
\vsp
\begin{equation*}
\begin{split}
  x &= (-1)^\sigma \ (d_0.d_1 d_2 d_3\cdots)_{\beta} \times \beta^e \\
  &=(-1)^\sigma (d_0+ d_1\beta^{-1}+d_2\beta^{-2}+\cdots)\times \beta^e
\end{split}
\end{equation*}
where $\sigma\in\{0,1\}$, {\em digits} $d_k\in\{0,1,\ldots,\beta-1\}$, and {\em exponent} $e$ is an integer. In computers we truncate the fraction and approximate $x$ by 
\vsp
\begin{equation}\label{flRep:eq}
  \wh x = (-1)^\sigma \ (d_0.d_1 d_2 \cdots d_{p-1})_{\beta} \times \beta^e 
\end{equation}
with $d_0\neq 0$ (for normalization) and exponent $e$ limited by 
$$
L\leqslant e\leqslant U.
$$
The string $d_0d_1\cdots d_{p-1}$ is called {\em mantissa} or {\em significand} and the portion $d_1\cdots d_{p-1}$ of the mantissa is called the {\em fraction}. Also, $p$ is called the {\em precision} of the floating-point system.
The key fact here is that in the floating point format the place of the base-$\beta$ point can be
changed by adjusting the exponent.
The condition $d_0\neq0$ makes the representation unique. For example, among all of the following representations of the decimal number
$123.4$, e.g. \vsp
  \[
  (1.234)_{10} \times 10^{2}  =   (0.1234)_{10} \times 10^{3} =  (0.01234)_{10} \times 10^{4},
  \]
the representation $(1.234)_{10} \times 10^{2}$ is permitted. 
The set of all real numbers that can be expressed in the form \eqref{flRep:eq}, with $d_0 \neq 0$, is denoted by $\F(\beta, p, L, U)$ and is referred to as the set of \textbf{normalized} floating-point (or machine) numbers. When the parameters are either fixed, known, or irrelevant to the context, this set is simply denoted by $\F$.
Numbers in
$\F$ are symmetric with respect to zero. The limited range of the exponent implies that $\wh x$ is limited
in magnitude to an interval which is called the {\em range} of the floating-point system. If $\wh x$ is
larger in magnitude than the largest number in the set $\F$, then $\wh x$ cannot be represented at all.
The same is true, in a sense, of numbers smaller than the smallest nonzero number in $\F$.
\vsp 

\begin{example} \label{ex:toy}
We generate the full list of all normalized numbers in the set
$\F(2, 3,-2, 1)$ and express each number in base $10$.
Such numbers are represented as $\pm(d_0.d_1d_2)_2\cdot 2^e$ for $e=-2,-1,0,1$.
Starting from the smallest exponent value $e = -2$, we obtain the following members:
\begin{align*}
\pm(1.00)_2 \times 2^{-2} &= \pm(1 \times 2^0 + 0 \times 2^{-1} + 0 \times 2^{-2}) \times 2^{-1} = \pm \frac{4}{16}= \pm0.25 ,\\
\pm(1.01)_2 \times 2^{-2} &= \pm(1 \times 2^0 + 0 \times 2^{-1} + 1 \times 2^{-2}) \times 2^{-1} = \pm\frac{5}{16} = \pm0.3125,\\
\pm(1.10)_2 \times 2^{-2} &= \pm(1 \times 2^0 + 1 \times 2^{-1} + 0 \times 2^{-2}) \times 2^{-1} = \pm\frac{6}{16} = \pm0.375 ,\\
\pm(1.11)_2 \times 2^{-2} &= \pm(1 \times 2^0 + 1 \times 2^{-1} + 1 \times 2^{-2}) \times 2^{-1} =  \pm \frac{7}{16} = \pm0.4375.
\end{align*}
Moving to cases $e = -1$, $e=0$ and $e=1$ we obtain numbers $\{\pm \frac{8}{16},\pm\frac{10}{16},\pm \frac{12}{16},\pm\frac{14}{16}\}$,
$\{\pm \frac{16}{16},\pm\frac{20}{16},\pm\frac{24}{16},\pm\frac{28}{16}\}$ and
$\{\pm \frac{32}{16},\pm\frac{40}{16},\pm\frac{48}{16},\pm\frac{56}{16}\}$, respectively. The positive side of this set is illustrated in Figure \ref{fig:toynumbers}.
As we can see, the smallest positive normalized number
is $(1.00)_2 \times 2^{-2} = 4/16 = 0.25$
and the largest number is $(1.11)_2 \times 2^{+1} = 56/16 = 3.5$.

\begin{center}
\begin{center}
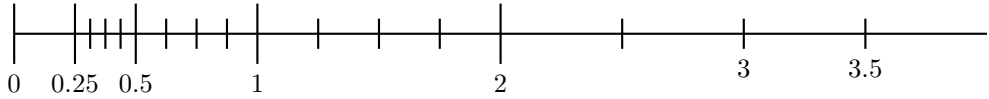

\begin{tikzpicture}
    \draw[color=black,thick] (0,0)--(12.8,0)  node[right] {};
    \draw[color=black,thick] (0,0.4)--(0,-0.4)  node[below] {\fs${0}$};
    \draw[color=black,thick] (0.8,0.4)--(0.8,-0.4)  node[below] {\fs${0.25}$};
    \draw[color=black,thick] (1.0,0.2)--(1.0,-0.2)  node[below] {};
    \draw[color=black,thick] (1.2,0.2)--(1.2,-0.2)  node[below] {};
    \draw[color=black,thick] (1.4,0.2)--(1.4,-0.2)  node[below] {};
    \draw[color=black,thick] (1.6,0.4)--(1.6,-0.4)  node[below] {\fs${0.5}$};
    \draw[color=black,thick] (2.0,0.2)--(2.0,-0.2)  node[below] {};
    \draw[color=black,thick] (2.4,0.2)--(2.4,-0.2)  node[below] {};
    \draw[color=black,thick] (2.8,0.2)--(2.8,-0.2)  node[below] {};
    \draw[color=black,thick] (3.2,0.4)--(3.2,-0.4)  node[below] {\fs${1}$};
    \draw[color=black,thick] (4.0,0.2)--(4.0,-0.2)  node[below] {};
    \draw[color=black,thick] (4.8,0.2)--(4.8,-0.2)  node[below] {};
    \draw[color=black,thick] (5.6,0.2)--(5.6,-0.2)  node[below] {};
    \draw[color=black,thick] (6.4,0.4)--(6.4,-0.4)  node[below] {\fs${2}$};
    \draw[color=black,thick] (8.0,0.2)--(8.0,-0.2)  node[below] {};
    \draw[color=black,thick] (9.6,0.2)--(9.6,-0.2)  node[below] {\fs${3}$};
    \draw[color=black,thick] (11.2,0.2)--(11.2,-0.2)node[below] {\fs $3.5$};
\end{tikzpicture}

\end{center}
\captionof{figure}{Positive normalized floating-point numbers in $\F(2,3,-2,1)$}\label{fig:toynumbers}
\end{center}
\end{example}
\vsp

 Floating-point numbers are not uniformly distributed throughout their range, but
are equally spaced only between successive powers of $\beta$. This means that
the distance between consecutive points is the same for a fixed exponent $e$ but increases as $e$ is increased.
The smallest and the largest positive numbers are denoted by $x_{min}$ and $x_{max}$, respectively.
\vsp 

\begin{workout}
Count the number of normalized machine numbers in $\F(\beta,p,L,U)$. Determine $x_{min}$ and $x_{max}$ in this set.
\end{workout}
\vsp

\begin{remark} \label{ulpnote2:rem}
For two arbitrary exponents $e$ and $e'$ with $L \leq e,e' \leq U$, the  the cardinality of both sets
$\F(\beta, p,L, U) \cap [\beta^e, \beta^{e+1}]$
and $\F(\beta, p,L, U) \cap [\beta^{e'}, \beta^{e'+1}]$ is the same. You can check this property in the toy example!
\end{remark}
\vsp 

The floating-point representation of number $1$ in base $\beta$ is $1 = (1.0 0 \cdots 0 0)_\beta \times \beta^0$. The smallest machine number larger than $1$ is equal to \vsp
\[
(1.0 0 \cdots 0 1)_\beta \times \beta^0 = 1 \times \beta^0 + 0 \times \beta^{-1} + \cdots + 0 \times \beta^{-(p-2)} + 1 \times \beta^{-(p-1)} = 1 + \beta^{-(p-1)}.
\]
The distance between these two successive floating-point numbers determines the precision of the machine and reflects the relative error of approximation in the system.
\vsp 
\begin{definition}
The {distance} between $1$ and the smallest machine number larger than $1$ is called
{\bf machine epsilon} and is denoted with $\varepsilon_{M}$.
The value of machine epsilon is
$
 \varepsilon_{M} = \beta^{-(p-1)}
$.
\end{definition}
\vsp 

The machine epsilon plays an important role in the analysis of rounding errors in numerical algorithms.
The distance between {\em any} machine number $x \in \F(\beta, p,L, U)$ and its consecutive machine number is called \textbf{unit in the last place} and is denoted by $\mathrm{ulp}(x)$. 
It states the weight of the last digit in the mantissa of the normalized number $x$.
Without loss of generality, we only consider positive numbers; an analogous result holds for negative numbers. 
If $x$ has representation 
$x = (d_0.d_1 d_2 \cdots d_{p-1})_{\beta} \times \beta^e$ and if all its digits are not $\beta-1$, then the next number, say $x^+$, is
$$
x^+ = (d_0.d_1 d_2 \cdots d_{p-1}+0.00\cdots01)_{\beta} \times \beta^e= x + \beta^{-(p-1)}\times 
\beta^e.
$$
This shows that 
\[
\mathrm{ulp}(x) = \beta^{e-p+1} = \varepsilon_{M} \beta^e.
\vsp
\]
If all digits of $x$ are $\beta-1$ then $x^+ = (1.00\cdots 0)_\beta \times \beta^{e+1}$. Now, $x^+$ is the first number in the next  interval (with a new exponent; for example floating-point numbers $0.5$, $1$ and $2$ in the toy example). In this case, we again can show that $\mathrm{ulp}(x) = \beta^{e-p+1}$.


Note that if $x > 0$, then $\mathrm{ulp}(x)$ is equal to the distance between $x$ and the smallest machine number larger than $x$ while if $x < 0$, then $\mathrm{ulp}(x)$ is the distance between
$x$ and the largest machine number smaller than $x$.
The concept of ulp is just a generalization of machine epsilon as $\varepsilon_{M} = \mathrm{ulp(1)}$.
\vsp

\begin{example}
In the toy example since the machine number $2.5$ belongs to $[2^1, 2^2]$ which corresponds to exponent $e = 1$ we have
$
 \mathrm{ulp(2.5)} = 2^{1 - 3 + 1}  = 0.5
$.
Similarly, the exponent of the machine number $x = 1.75$ is $0$ thus
$
 \mathrm{ulp(1.75)} = 2^{0 - 3 + 1}  = 0.25
$.
Similarly, ulp$(1)=\varepsilon_M = 0.25$.
\end{example}
\vsp 

\begin{remark}
 Since computers use base $2$ with digits $\{0,1\}$, the only choice for the leading bit $d_0$ in the normalized binary floating-point numbers is
 $d_0 = 1$. This creates an opportunity to save some memory if we avoid occupying the bit $d_0$ with the fixed value of $1$ by imposing it implicitly. 
This implicit bit whose value is always $1$ for any normalized binary number is called the \textbf{hidden bit}.
 \end{remark}
\vsp 

\subsection{Rounding modes} \label{rounding_modes:subsec}

In computers, a real number $x$ must be mapped in to (approximated by) a `nearby' machine number before we can start any computation with it. More precisely, we need a map \vsp
     \[
      \fl: \mathbb{R} \to \F(\beta, p,L, U)
      \]
from the uncountable set of real numbers to the set of machine numbers. Such mappings are called {\em rounding modes}.
The process of choosing a nearby floating-point number $\fl(x)$ to approximate
a given real number $x$ is called {\em rounding}, and the error introduced by such
an approximation is called {\em rounding error}, or {\em roundoff error}. Two commonly used
rounding rules are {\bf chopping} and {\bf rounding to nearest}\footnote{Two additional rounding modes are {\em rounding to up} (or round toward $\infty$) and {\em rounding to down} (or round toward $-\infty$) which are used in interval arithmetic. In interval arithmetic, numbers are represented by intervals (an upper and lower bound) rather than a single value. For example, $\pi$ may be approximated by [3.141, 3.142].}.
Assume that $x = \pm (d_0.d_1\cdots d_{p-1}d_p\cdots)_\beta\times \beta^e$.
In chopping rule, the mantissa of $x$ is truncated after the $(p - 1)-$st digit to get \vsp
$$
\fl(x) = \pm (d_0.d_1\cdots d_{p-1})_\beta\times \beta^e.
$$
Since
$\fl(x)$ is indeed the next floating-point number towards zero from $x$, this rule is 
sometimes called {\em round toward zero}. In rounding to nearest we have \vsp
$$
\fl(x) = \pm (d_0.d_1\cdots d_{p-2}\tilde d_{p-1})_\beta\times \beta^e,
$$
where $\tilde d_{p-1}$ is either $d_{p-1}$ or $d_{p-1}+1$ depending on to which one of the consecutive floating-point numbers $x$ is closer. If $x$ falls exactly midway then we have a {\em tie}. In this case $x$ is rounded to the nearest floating-point number with an even last significant digit. This rule is called {\em round to nearest ties to even}.

Clearly, $\fl(x) = x$
if $x \in \F$. Moreover, $\fl(x) \leqslant \fl(y)$ if $x \leqslant y$ for all $x, y \in \R$ (monotonicity
property). Since the distance between two consecutive machine numbers is $\varepsilon_M\beta^{e}$, the relative error of chopping can be bounded as 
\vsp
$$
\frac{|\fl(x)-x|}{|x|} \leqslant \frac{\varepsilon_M\beta^{e}}{(d_0.d_1d_2\cdots)_\beta\beta^{e}}\leqslant \frac{\varepsilon_M}{(1.0)_\beta} = \varepsilon_M.
\vsp
$$
Similarly, for rounding to nearest we have 
\vsp 
$$
\frac{|\fl(x)-x|}{|x|} \leqslant \frac{\varepsilon_M}{2}.
\vsp
$$
Note that, in a floating-point system both large and small numbers are represented with
nearly the same relative precision.
The relative roundoff error for rounding to nearest is half of that for chopping, the main reason rounding to nearest is
the default rounding rule in all standard systems although it is more
expensive to implement. From here on, by `rounding' we mean `rounding to nearest'.
\vsp 

\begin{theorem}
In the floating-point number system $\F(\beta, p, L, U)$ every real number in the
floating-point range can be represented with a relative error, which does not exceed
the {\bf unit roundoff} $u$, which is defined by
$$
u = \begin{cases} {\varepsilon_M}/{2}, & \mathrm{if~rounding~is~used}\\ \varepsilon_M, & \mathrm{if~chopping~is~used} \end{cases}
$$
\end{theorem}
\vsp

The quantity $u$ is a natural unit for
relative changes and relative errors. For example, termination criteria in iterative methods
usually depend on the unit roundoff. The following theorem is the fundamental theorem of rounding errors.
\vsp 

\begin{theorem}
If $x\in\R$ is such that $|x|\in [x_{min},x_{max}]$, then
\begin{equation}\label{fl_rule}
\fl(x) = x(1+\delta), \quad |\delta|\leqslant u.
\end{equation}
\end{theorem}
\vsp 
The proof of this theorem follows the error bounds on chopping and rounding to nearest mappings. 
The fundamental rule \eqref{fl_rule} will be used to analyze the effect of initial and intermediate rounding errors in the final result of an algorithm.
\vsp 

\subsection{Subnormal numbers}

Recall the toy example once more. Due to the normalization requirement 
$ d_0 \neq 0$, there is a relatively large gap between zero and the smallest positive number. Additionally, zero cannot be represented by the normalized floating-point representation \eqref{flRep:eq}. These have
an unfortunate impact on the validity of some of the most important algebraic properties when performing arithmetic with machine numbers in $\F$.
For instance, in the current scenario, there could be two different normalized machine numbers, $x$ and $y$, such that $x - y$ falls within this gap and is consequently approximated by zero!
A remedy for this situation is to allow the leading digit $d_0$ to be zero but
only when the exponent is at its minimum value $L$. Then the gap around zero can
be filled in by additional floating-point numbers which are called {\bf subnormal} or denormalized numbers.
A subnormal floating-point number is of the form 
\vsp
\begin{equation*}
  \wh x = \pm (0.d_1d_2\cdots d_{p-1})_\beta \beta^{L}, \quad d_k\in\{0,1,\ldots,\beta-1\}.
\end{equation*}
The fixed exponent $e=L$ makes this representation unique.
\vsp 

\begin{example}
For the toy example subnormal numbers are
\vsp 
$$
\pm(0.00)_22^{-2} = \pm0,\quad \pm(0.01)_22^{-2} = \pm\frac{1}{16},\quad\pm(0.10)_22^{-2} = \pm\frac{2}{16},\quad\pm(0.11)_22^{-2} = \pm\frac{3}{16}.
$$
The non-negative normalized and denormalized numbers are depicted in Figure \ref{fig:subnormals}.
\begin{center}
\begin{center}
\begin{tikzpicture}
    \draw[color=black,thick] (0,0)--(12.8,0)  node[right] {};
    \draw[color=red,thick] (0,0.4)--(0,-0.4)  node[below] {\fs${0}$};
    \draw[color=red,thick] (0.2,0.2)--(0.2,-0.2)  node[below] {};
    \draw[color=red,thick] (0.4,0.2)--(0.4,-0.2)  node[below] {};
    \draw[color=red,thick] (0.6,0.2)--(0.6,-0.2)  node[below] {};
    \draw[color=black,thick] (0.8,0.4)--(0.8,-0.4)  node[below] {\fs${0.25}$};
    \draw[color=black,thick] (1.0,0.2)--(1.0,-0.2)  node[below] {};
    \draw[color=black,thick] (1.2,0.2)--(1.2,-0.2)  node[below] {};
    \draw[color=black,thick] (1.4,0.2)--(1.4,-0.2)  node[below] {};
    \draw[color=black,thick] (1.6,0.4)--(1.6,-0.4)  node[below] {\fs${0.5}$};
    \draw[color=black,thick] (2.0,0.2)--(2.0,-0.2)  node[below] {};
    \draw[color=black,thick] (2.4,0.2)--(2.4,-0.2)  node[below] {};
    \draw[color=black,thick] (2.8,0.2)--(2.8,-0.2)  node[below] {};
    \draw[color=black,thick] (3.2,0.4)--(3.2,-0.4)  node[below] {\fs${1}$};
    \draw[color=black,thick] (4.0,0.2)--(4.0,-0.2)  node[below] {};
    \draw[color=black,thick] (4.8,0.2)--(4.8,-0.2)  node[below] {};
    \draw[color=black,thick] (5.6,0.2)--(5.6,-0.2)  node[below] {};
    \draw[color=black,thick] (6.4,0.4)--(6.4,-0.4)  node[below] {\fs${2}$};
    \draw[color=black,thick] (8.0,0.2)--(8.0,-0.2)  node[below] {};
    \draw[color=black,thick] (9.6,0.2)--(9.6,-0.2)  node[below] {\fs${3}$};
    \draw[color=black,thick] (11.2,0.2)--(11.2,-0.2)node[below] {\fs $3.5$};
\end{tikzpicture}
\end{center}

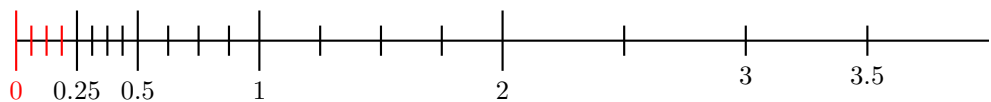
\captionof{figure}{Non-negative normal and subnormal floating-point numbers in $\F(2,3,-2,1)$}\label{fig:subnormals}
\end{center}
\end{example}
\vsp
\begin{workout}
Count the number of subnormal floating-point numbers in $\F(\beta,p,L,U)$. Determine the distance between two successive subnormal numbers.
\end{workout}
\vsp 

The introduction of subnormal numbers guarantees that the subtraction of (nearby) floating-point numbers (with the same sign or the addition of floating-point numbers with opposite signs) never gives zero, i.e., the essential relation
\[
x  = y \iff x - y = 0
\]
is now valid for $x, y \in \F$. This is caused by the fact that subnormal numbers can represent the non-zero distance between two nearby floating-point numbers. Several examples of how subnormal numbers make
writing reliable floating-point code easier are analyzed in \cite{Demmel:1984}.
\vsp

\begin{remark}
Even with subnormal numbers, the validity of the formula
\[
x  = y \iff x - y = 0
\]
is not guaranteed if $x, y \notin \F$. 
\end{remark}
\vsp
Note that, subnormal numbers have inherently lower precision than
normalized numbers because they have fewer significant digits in their fractional
parts.
\vsp 
\begin{remark}
If $x$ is the result of an operation on numbers of $\F$ and  $x \in (-\infty,-x_{max}) \cup (x_{max},\infty)$
then $\fl(x)$ can not be defined. 
On the other side, if $x \in (-x_{min}, x_{min})$, the operation of rounding is defined anyway
(even in absence of subnormal numbers). The first case is referred to {\bf overflow} and the
second case to {\bf underflow}. The values of $x_{min}$ and $x_{max}$ are called the underflow level (UFL) and the overflow level (OFL), respectively.

In the presence of subnormal numbers, we can think of {\bf gradual underflow} instead of underflow itself. Because
if $x$ is approximated by a subnormal number with representation
$\fl(x)=\pm(0.0\cdots 0d_{p-k}\cdots d_{p-1})_\beta \beta^L$ with $d_{p-k}\neq 0$ for an integer $k$ with $1\leqslant k< p$, then
\\
\begin{equation*}
  \frac{|\fl(x)-x|}{|x|}\leqslant \frac{1}{2}\frac{\beta^{L-(p-1)}}{(0.0\cdots 010\cdots 0)_\beta \beta^L}=\frac{1}{2}\frac{\beta^{-(p-1)}}{\beta^{-(p-k)}} = \frac{1}{2}\beta^{-k+1}.
\end{equation*}
\\
This shows that the relative roundoff error increases {\em gradually} as $k$ (the number of significant digits of $x$) decreases, or equivalently as $\fl(x)$ approaches zero.
Remember that in the absence of subnormal numbers we have underflowing to zero as soon as $x$ falls into the underflow region.
\end{remark}
\vsp

\subsection{IEEE standard for floating point arithmetic}
In 1941, Konrad Zuse built the first computer in Berlin, which was called Z3. It was an electro-mechanical computer with a binary ($\beta = 2$) system for machine numbers occupying a total of 22 bits. Z3 was destroyed in a bombardment of Berlin during the World War II. It is often said that the first fully electronic computer was ENIAC\footnote{Electronic Numerical Integrator and Computer}, built between 1943 and 1945 at the University of Pennsylvania\footnote{January 15 is known as ENIAC Day and is celebrated annually in the United States.}. Subsequently, more efficient and smaller computers were developed.
During 1960-1980, floating point computation was used as a basis of scientific computation. However, each computer manufacturer developed its own floating point system. See Table \ref{tab:computers} for a few examples.
\vsp 
\begin{table}[!h]
\begin{center}
\caption{Different formats on different computers}\label{tab:computers}
 \begin{tabular}{|l c c c|}
 \hline
 computer & $\beta$ & $p$ & $U = -L$  \\ [0.5ex]
 \hline
 $\rm{IBM\ 7090}$  &$2$ & $27$ & $2^7$  \\
 $\rm{Borroughs\ 5000\ Series}$ & $8$ & $13$ &$2^6$ \\
 $\rm{IBM\ 360/370}$ & $16$ &  $6$ & $2^6$\\
 $\rm{DEC\ 11/780\ VAX}$ & $2$ & $24$ & $2^7$ \\
$\rm{Hewlett\ Packard\ 67}$& $10$ & $10$ & $99$ \\ [1ex]
 \hline
\end{tabular}
\end{center}
\end{table}

Due to differences between various implementations of floating-point systems, computer programs were hardly portable and often yielded  different results when running the same code on different machines. A program that produced accurate results on one computer might not even run on another. These discrepancies were caused not only by variations in parameters such as \(\beta\), \(p\), \(L\), and \(U\), but also by differences in the implementation details of floating-point arithmetic.
To address these issues, in a great cooperation between academic computer scientists and hardware designers, a standard for binary floating point representation and arithmetic was developed in 1985 which was supported by the Institute of Electrical and Electronics Engineers
(shortly, IEEE). The standard was called {\em IEEE 754}. In 1985, a second standard was established by IEEE (called IEEE 854) for both decimal and binary floating-point systems\footnote{The leader of the academic computer scientists who developed the standard was William Kahan from the University of California, Berkeley. He received the ACM (Association for Computing Machinery) Turing Award in 1989 for his contributions to this standard.}.

\begin{figure}[!th]
\centering
\includegraphics[scale=0.6]{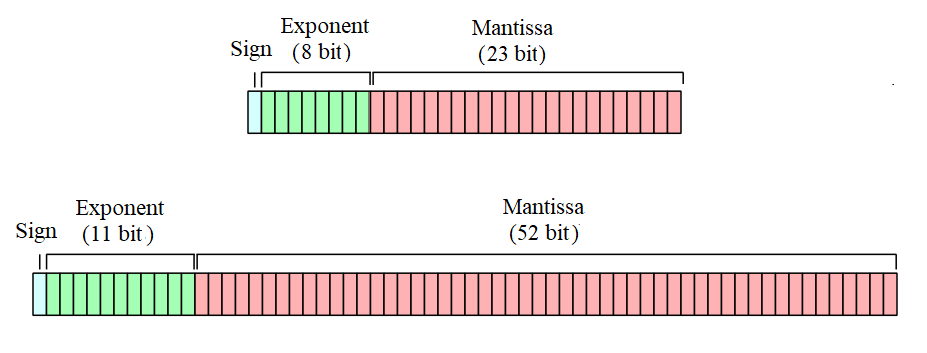}
\caption{Number of bits in single (top) and double (down) precisions in the IEEE standard.}
\label{fig:ieee_bits}
\end{figure}

The IEEE standard specifies two basic representation formats, {\em single} or fp32, and {\em double} or fp64.
The general structure of the two data types single and double are the same and the only difference is in the number of bits used to represent the mantissa and the exponent. See Figure \ref{fig:ieee_bits}. 

In the single data type (fp32), there are a total of 32 bits, the first one of which is called the sign bit and represents the sign of the number, the next 8 bits correspond with the exponent of the number, and the remaining 23 bits are used to represent the mantissa. Following the idea of the {\em hidden bit}, the precision of the single format is therefore equal to $p = 23 + 1$. These numbers are increased to 11 bits for exponent, 52 bits for mantissa with precision $p=52+1$ for the double format.
A {\em biased exponent} is stored and
no sign bit used for the exponent which prevent wasting one bit to represent the sign of the exponent. In single precision, $L=-126$ and $U = 127$, and $e + 127$ is stored in eight bits. This range includes 254 integers, which is two fewer than the total of $2^8 = 256$ numbers that can be represented with eight bits. The other two exponents, $-127$ and $+128$, are reserved for special numbers. 
In double precision,
$L=-1022$, $U=1023$ and the biased exponent $e+1023$ is stored in $11$ bits. Again two extra representable exponents $-1023$ and $1024$ are reserved for storing special numbers. 

The IEEE standard includes {\em extended single} and {\em extended double} formats
that offer extra precision and exponent range.
The characteristics of all formats are summarized in Table \ref{tab:fp_formats}.

\begin{table}[!h]
\begin{center}
\caption{Single and double and extended precision IEEE formats}\label{tab:fp_formats}
\begin{tabular}{|lrrrrr|}
  \hline
  Format & Bits No. & $p$~ & $e$~~~ & $L$~~~~ & $U$~~~~ \\
  \hline
  Single & $32$~bits & $24$ & $8$~bits & $-126$ & $127$ \\
  Extended single & $43$~bits & $32$ & $11$~bits & $ -1022$ & $ 1023$ \\
  Double  & $64$~bits & $53$ & $11$~bits & $-1022$ & $1023$ \\
  Extended double & $79$~bits & $64$ & $15$~bits & $ -16,382$ & $ 16,383$ \\
  \hline
\end{tabular}
\end{center}
\end{table}
As the double precision satisfies the
requirements for the extended single format, so three precisions single, double and extended double suffice.
Extended double format can for instance be used in intermediate calculations for computing of elementary functions accurately in the double precision.

In IEEE standard, an exponent $e=L-1$ (for example $e=-127$ in single precision) and a nonzero mantissa corresponds to a subnormal number $\pm(0.d_1\cdots d_{p-1})_22^{L}$.
There are distinct representations for $+0$ and $-0$ with
exponent $e=L-1$, and a zero mantissa\footnote{
The signed zero distinguishes between positive and negative underflowed numbers.
Another use of singed zero is in the computation of complex elementary functions.}. This means that the representation of $+0$ in the single precision is 
\vsp
$$
0|00000000|00000000000000000000000
$$
Note that the biased exponent $e+127 = -127+127=0$ is stored in eight exponent bits.

Infinity (overflow) is also signed and $\pm\infty$ is represented by the exponent $e=U+1$ and a zero
mantissa. It is also obtained from operation $x/0$ for $x\neq0$. The IEEE infinity obeys the mathematical rules
$\infty+\infty=\infty$, $(-1)\times \infty=-\infty$ and $x/\infty = 0$.

The standard also introduces special numbers called \verb+NaN+ (Not a Number) for undefined floating point operations like $0/0$, $\infty-\infty$ and $0\times\infty$. A \verb+NaN+ is stored with exponent $U+1$ and a nonzero mantissa. When a \verb+NaN+ and an ordinary floating-point number are combined the result is another \verb+NaN+. A \verb+NaN+ is also often used for uninitialized or missing data.
Table \ref{tab:IEEE_rep} summarizes all possible representations in IEEE 754 standard.

\begin{table}[!h]
\begin{center}
\caption{IEEE 754 representation. Here $m=d_1\cdots d_{p-1}$ is mantissa}
\label{tab:IEEE_rep}
\begin{tabular}{|ccll|}
  \hline
Type of Number &Exponent & Mantissa & Represents\\
  \hline
Zero &$e=L-1$ & $m=0$ & $\pm0$\\
Subnormal  &$e=L-1$ & $m\neq0$ & $\pm(0.m)_22^L$\\
Normal & $L<e<U$ & $m\neq 0$ or $m=0$ & $\pm(1.m)_22^e$\\
Infinity & $e=U+1$& $m=0$ & $\pm\infty$\\
Not a Number& $e=U+1$& $m\neq0$& \verb+NaN+\\
\hline
\end{tabular}
\end{center}
\end{table}
\begin{example}
To represent the decimal number \(89.75\) in IEEE standard single format, we start by converting it to the binary system. We have \(89 = (1011001)_2\) and \(0.75 = (0.11)_2\). This gives
\[
89.75 = (1011001.11)_2 = (1.01100111)_2 \times 2^6.
\]

This number is a positive normalized floating-point number (with the hidden bit \(d_0 = 1\)), an exponent \(e = 6\), and a mantissa \(m = 01100111\). 
We store the biased exponent \(e + 127 = 133\), which in binary is \( (10000101)_2\). Thus, the final IEEE 754 representation is
\vsp
\[
0|10000101|01100111000000000000000
\]

Note that the 14 remaining bits in the mantissa are filled with zeros.
\end{example}
\vsp 
\begin{workout}
Which decimal numbers do these two floating-point single precision arrays represent?
\begin{itemize}
\item[(a)] $1|01011001|01110100000000000000000$

\item[(b)] $1|00000000|01110100000000000000000$
\end{itemize}
\end{workout}
\vsp 

The machine epsilon, the smallest positive normalized number and the largest normalized number in the single precision are
\begin{align*}
  \varepsilon_M & = 2^{-(p-1)} = 2^{-23} \approx 1.19 \times 10^{-7}\\
  x_{min} & = (1.00\cdots 0)_22^{-126}\approx 1.18\times 10^{-38}\\
  x_{max} & = (1.11\cdots 1)_22^{+127}\approx  3.40 \times 10^{+38}.
\end{align*}
These numbers are refined to
\begin{align*}
  \varepsilon_M & = 2^{-52} \approx 2.22 \times 10^{-16}\\
  x_{min} & = (1.00\cdots 0)_22^{-1022}\approx 2.23\times 10^{-308}\\
  x_{max} & = (1.11\cdots 1)_22^{+1023}\approx  1.80 \times 10^{+308}
\end{align*}
in the double precision format.
The double data type (fp64) is the default format in several software for scientific computing.

In Python the following commands print all machine parameters in double precision. 
\begin{shaded}
\vspace*{-5mm}
\begin{verbatim}
import numpy as np
print(np.finfo(float))
\end{verbatim}
\vspace*{-5mm}
\end{shaded}
\noindent
The output is:
\begin{shaded}
\vspace*{-4mm}
\begin{verbatim}
Machine parameters for float64
---------------------------------------------------------------
precision =  15   resolution = 1.0000000000000001e-15
machep =    -52   eps =        2.2204460492503131e-16
negep =     -53   epsneg =     1.1102230246251565e-16
minexp =  -1022   tiny =       2.2250738585072014e-308
maxexp =   1024   max =        1.7976931348623157e+308
nexp =       11   min =        -max
---------------------------------------------------------------
\end{verbatim}
\vspace*{-5mm}
\end{shaded}
\noindent
Here, \verb+eps+ stands for machine epsilon, \verb+epsneq+ for unit roundoff, \verb+tiny+ for the smallest positive normalized numbers, \verb+max+ for the largest normalized number, \verb+nexp+ for number of allocated bits for the exponent, \verb+minexp+ for $L$ and \verb+maxexp+ for $U+1$.
The decimal \verb+precision+ $15$ is reported here instead of the binary precision $53$.
To access the variables, for example to access the unit roundoff, we can write
\begin{shaded}
\vspace*{-5mm}
\begin{verbatim}
u = np.finfo(float).epsneg
\end{verbatim}
\vspace*{-5mm}
\end{shaded}
\noindent
To observe the corresponding parameters for the single precision format write
\begin{shaded}
\vspace*{-5mm}
\begin{verbatim}
print(np.finfo(np.float32))
\end{verbatim}
\vspace*{-5mm}
\end{shaded}
\noindent
Note that the comments above are provided in \verb+numpy+ Python library as we imported this library with abbreviation \verb+np+ in the first input line. The ulp of a number can be obtained by command \verb+spacing+ from \verb+numpy+ module or the command \verb+ulp+ from \verb+math+ module.
\begin{shaded}
\vspace*{-5mm}
\begin{verbatim}
In [1]: import numpy as np
In [2]: import math
In [3]: math.ulp(5)
Out[3]: 8.881784197001252e-16
In [4]: np.spacing(5)
Out[4]: 8.881784197001252e-16
In [5]: math.ulp(1)
Out[5]: 2.220446049250313e-16
In [6]: math.ulp(0)
Out[6]: 5e-324
\end{verbatim}
\vspace*{-5mm}
\end{shaded}
\noindent
As we observe, \verb+math.ulp(1)+ is another command to access the machine epsilon.
Also, \verb+math.ulp(0)+ gives the smallest positive subnormal representable floating-point number\footnote{See also \texttt{sys} module for these and other system-specific parameters and functions.}.
\vsp

\subsection{Disasters caused by inappropriate use of floating point arithmetic}

Even though the rounding errors in individual operations are typically small, their accumulation in complicated algorithms can lead to significant errors with potentially disastrous consequences. Here are a few historical events illustrating such cases.

{\bf Patriot missile failure in 1991 due to rounding errors\footnote{Source: Robert D. Skeel. Roundoff error and the patriot missile. SIAM News, 25(4): 11, Jul. 1992.}}:
On February 25, 1991, during the Iraq and Kuwait war, an American Patriot missile battery in Dharan, Saudi Arabia, failed to intercept an incoming Iraqi Scud missile. The Scud struck an American Army camp and killed 28 soldiers. The Patriot missile system was designed to track incoming missiles, predict their paths, and intercept them. However, it failed due to an inaccurate calculation of the time when the Patriot should have been launched.
The time was tracked by the system's internal clock in tenths of a second, which was then multiplied by \( 1/10 \) to convert it to seconds. This calculation was done using a 24-bit register. The value \( 1/10\), which has a non-terminating binary expansion, was rounded to $24$ bits after the radix point. Although the rounding error seemed negligible, it accumulated over time and led to a significant error.
By the time the Patriot battery had been operational for about $100$ hours, and an easy calculation shows that the resulting time error due to the rounding error was about $0.34$ seconds. The binary expansion of $1/10$ is
\begin{align*}
(0.1)_{10} &= 2^{-4} + 2^{-5} + 2^{-8} + 2^{-9} + 2^{-12} + 2^{-13} + \cdots \\
& = (1.10011001100 \cdots)_2 \times 2^{-4} = (1.100\overline{1100})_2 \times 2^{-4}.
\end{align*}
The $24$-bit register stored this value as \( 0.00011001100110011001100 \), introducing an error of approximately \( 0.000000095 \) in decimal. When multiplied by the number of tenths of a second in $100$ hours, this error became
$
0.000000095 \times 100 \times 60 \times 60 \times 10 = 0.34
$.
A Scud missile travels at about $1676$ meters per second, so travels more than half a kilometer in $0.34$ seconds. This error was sufficient to place the incoming Scud outside the ``range gate'' that the Patriot missile was tracking.

\begin{figure}[!th]
\centering
\includegraphics[scale=0.0705]{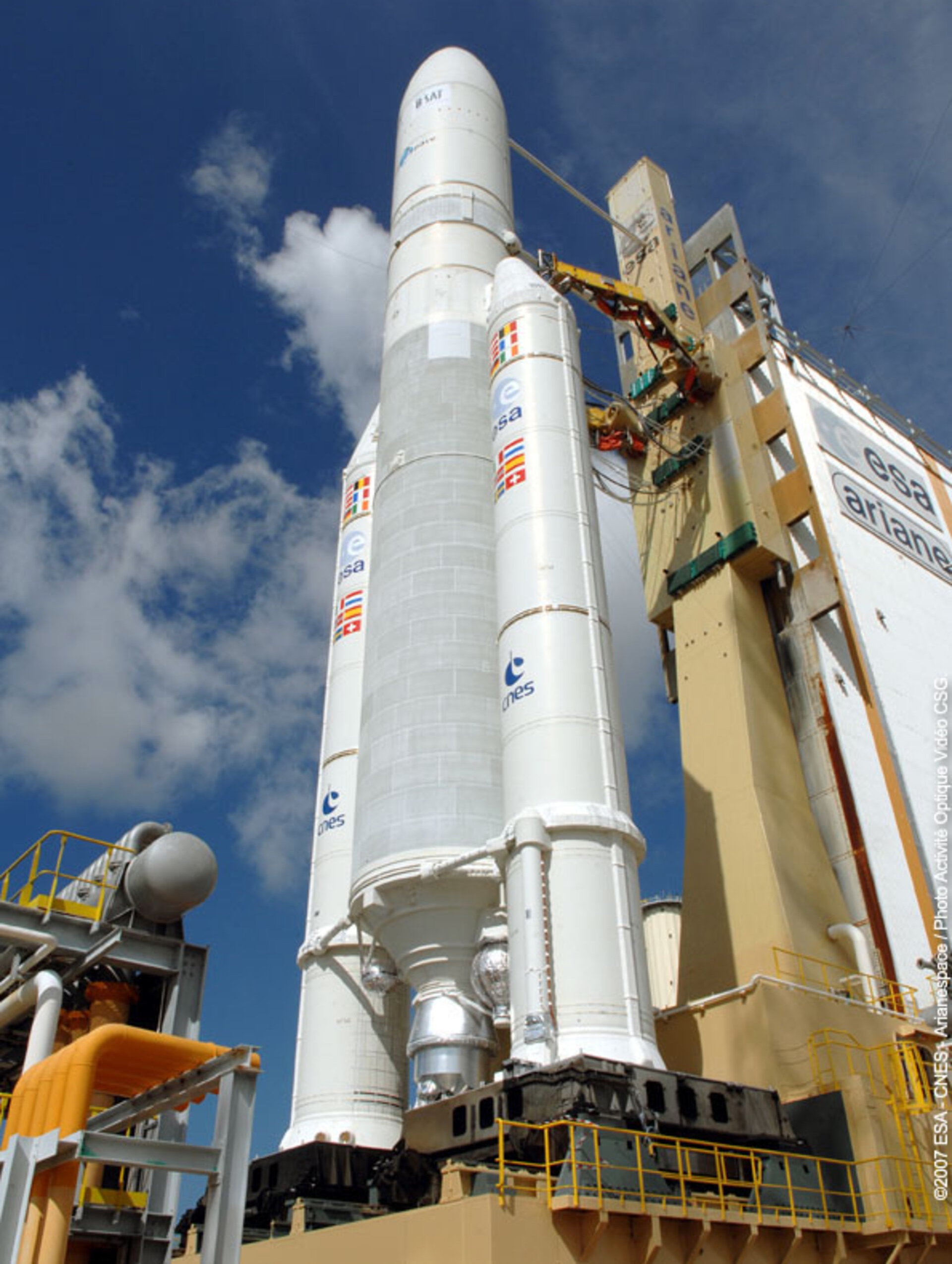} \includegraphics[scale=0.6]{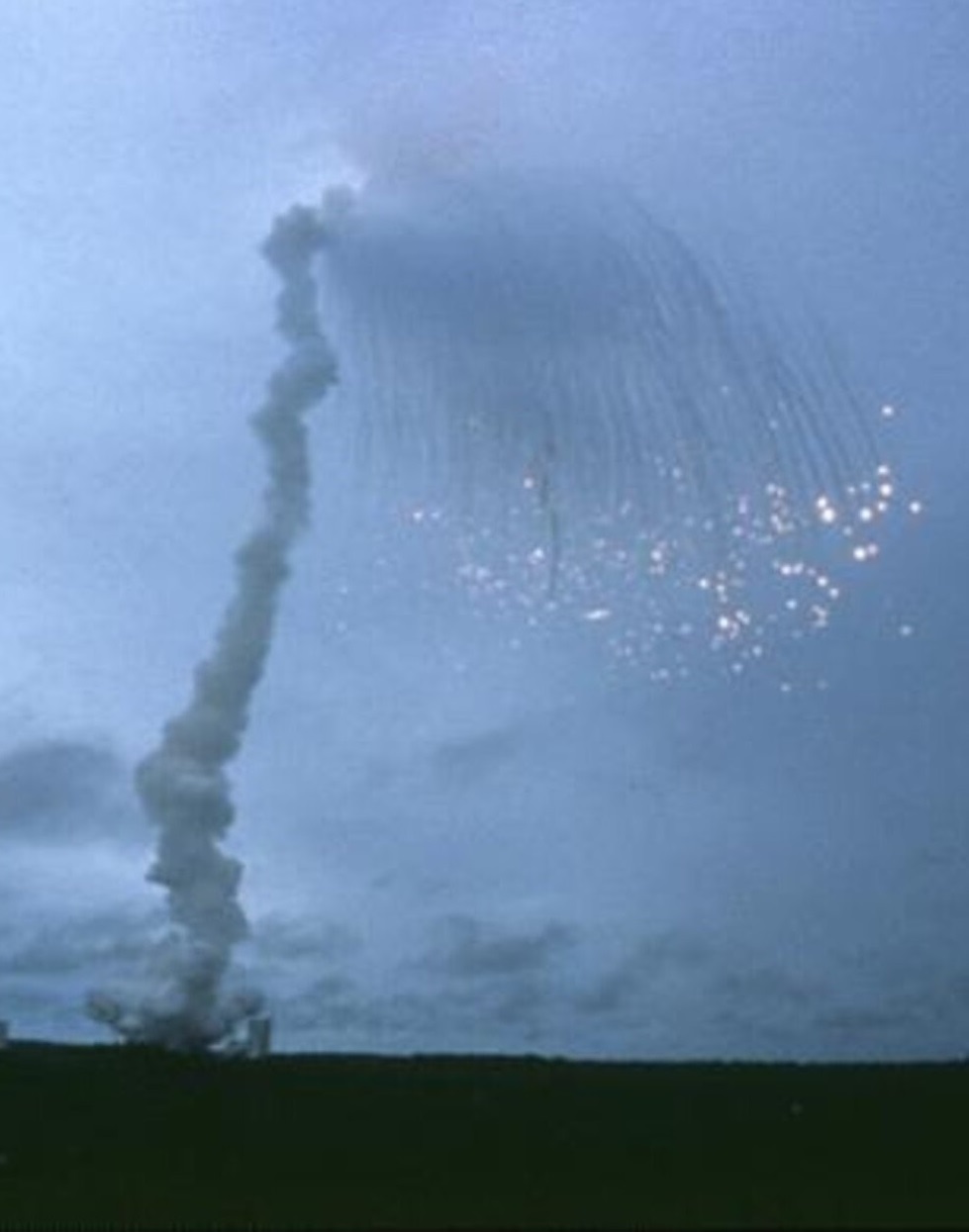}
\caption{Ariane 5 rocket (left), explosion after lift off (right). Photos from \texttt{https://www.esa.int}.}\label{fig:ariane}
\end{figure}

{\bf Explosion of the Ariane 5 rocket in 1996 due to an overflow error\footnote{Source: M. Dowson, The Ariane 5 Software Failure, ACM SIGSOFT Software Engineering Notes. 22 (1997): 84.}}:
On June 4, 1996, an Ariane 5 rocket launched by the European Space Agency exploded just 40 seconds after lift-off. This rocket was on its maiden voyage following a decade of development costing $\$$7 billion. The destroyed rocket and its cargo were valued at $\$$500 million.  It turned out that the cause of the failure was an overflow error in the inertial reference system. Specifically a 64 bit floating point number relating to the horizontal velocity of the rocket with respect to the platform was converted to a 16 bit signed integer. The number was larger than 32,768, the largest integer storeable in a 16 bit signed integer, and thus the conversion failed.
\vsp

\subsection{Algebraic properties of floating point arithmetic}
The set of real numbers is a {\em field} which makes the exact manipulation and analysis in $\R$ perfectly easy.
However, floating-point addition and multiplication are commutative but not associative, and the distributive law also fails for them. The same holds true for subtraction and division. This makes the analysis of floating-point computations quite difficult. Here we give
some examples to illustrate the situation.
\vsp 
\begin{example}
Let $x = 0.5, y= 2.5$ and $z=0.75$ in the toy example $\F(2,3,-2,1)$ with rounding to nearest ties to even. We then have
\begin{align*}
& \fl(0.5 + \fl(1.5 + 0.75)) = \fl(0.5 + \fl(2.25)) = \fl(0.5 + 2) = \fl(2.5) = 2.5,\\
& \fl(\fl(0.5 + 1.5) + 0.75) =\fl( \fl(2) + 0.75) = \fl(2 + 0.75) = \fl(2.75) = 3,
\end{align*}
which shows that the floating-point addition is not associative.
\end{example}
\vsp 
\begin{example}
Let $x = 2, y = \varepsilon_M$, and $z = -2$ in IEEE double precision. Let us compare $(x+y)+z$ and $x+(y+z)$ by doing the operations in  Python:
\begin{shaded}
\vspace*{-0.5cm}
\begin{verbatim}
In [1]: import math
In [2]: eps = math.ulp(1)
In [3]: (2 + eps) - 2
Out[3]: 0.0
In [4]: 2 + (eps - 2)
Out[4]: 2.220446049250313e-16
\end{verbatim}
\vspace*{-0.5cm}
\end{shaded}
\noindent
Following the definition of the machine epsilon and because the distance between any two consecutive machine numbers in interval $[1,2]$ is the same, we conclude that the distance between any two consecutive machine numbers in $[1,2]$ is \texttt{eps}.
In addition, the distance between consecutive machine numbers $[2,4]$ is 2\texttt{eps}, i.e., the smallest machine number larger than 2 is \texttt{2+2eps}. Since
\texttt{2+eps} is not a machine number, we have to round it to the nearest machine number. \texttt{2+eps} is right in the middle of $2$ and \texttt{2+2eps} and since $2$ is an even number (the last bit of its mantissa in base two is zero) and so \texttt{2+2eps} is an odd number, we understand why we should get \texttt{(2+eps) = 2}. That is basically why we got zero as the result of \texttt{(2+eps) - 2}.
Furthermore, because the set of machine numbers is symmetric with respect to zero, the distance between consecutive machine numbers from interval $[-2,-1]$ is \texttt{eps}. Therefore,  the exact value of
\texttt{eps-2} is somewhere in the interval $[-2,-1]$ where the gap between consecutive machine numbers is \texttt{eps}. This means that \texttt{eps-2} is actually a machine number! So, there is no rounding errors in computing $y+z$. That is why we got \texttt{eps} when we add $x=2$ to $(y+z)$.
\end{example}
\vsp
\begin{workout}
Show with some examples that (a) floating point multiplication is not necessarily associative, (b)
floating point multiplication is not necessarily distributive over floating point addition, i.e. there exist
$x, y, z \in \F$ such that $\fl (x  \fl(y + z)) \neq \fl(\fl(xy) + \fl(xz))$.
Also, (c) there exist
$x, y, z \in \F$ such that $\fl(x+y) = \fl(x+z)$ while $y\neq z$, and (d) there exist
$x, y, z \in \F$ such that $\fl(xy) = \fl(xz)$ while $y\neq z$. Property (c) for $y=0$ and property (d) for $y=1$ are valid showing that addition and multiplication neutral elements are not unique in floating-point arithmetic.
\end{workout}
\vsp 

These all show that the order of floating-point operations could affect the accuracy of the result. We will see some examples in analysis of algorithms in section \ref{sect-analysis_algorithms}.
\vsp 

\subsection{Floating-point arithmetic models}
The set of real numbers is closed under the four basic arithmetic operations while the set of machine numbers is not. More precisely, if
$x,y \in \F$ and $\ast \in \{ +, -, \times, / \}$, then it may happen that
$x \ast y \notin \F$.
For example, consider the operation of dividing $1$ by $3$ on a machine whose representation base $\beta$ is either $2$ or $10$. We know that the exact value of $1/3$ has infinitely many digits in both binary and decimal representations. Thus the results need to be rounded to a floating-point number. The rounding rule \eqref{fl_rule} suggests the model
\vsp 
\begin{equation}\label{fl_op_rule}
  \fl(x\ast y) = (x\ast y)(1+\delta), \quad |\delta|\leqslant u,\quad \ast\in\{+,-,\times,/ \}
\end{equation}
for normalized numbers $x,y\in\F$ provided that $x\ast y$ is in the normalized range.
Sometimes the floating-point computation is more precise than what the
model \eqref{fl_op_rule} assumes. An obvious example is that when $x\ast y\in \F$, there is no rounding error at all.  
Among elementary functions, the roundoff error introduced in computing $\sqrt x$ in IEEE standard obeys the same model
\vsp 
\begin{equation*}
  \fl(\sqrt x) = \sqrt x(1+\delta), \quad |\delta|\leqslant u,\quad x\in\F.
\end{equation*}
It is shown in \cite{Higham:2002} that the model \eqref{fl_op_rule} is valid for addition and subtraction operators if the arithmetic is supported by the {\bf guard digit} for subtraction. Otherwise a weaker model holds true. The role of guard digit can be easily explained by the following simple example from the mentioned book. Consider a floating-point arithmetic system with base $\beta=2$ and precision $p=3$. Let us subtract from $1$ the next smaller floating-point number:
\vsp
$$
\begin{array}{rlcrll}
 &\;\;\;1.00 \times 2^0    &                &&\;\;\;1.00~ \times 2^0&\\
 &-1.11 \times 2^{-1} &\Longrightarrow & &-0.111 \times 2^0 &\\
\cline{4-5}
 &                   &                & &\;\;\;0.001 \times 2^0 &= 1.00\times 2^{-3}
\end{array}
\vsp
$$
Since for subtraction we have to scale both numbers to the same exponent (the larger one), a third digit is introduced in the mantissa of the second number. This digit, if not dropped, is known as {\em guard digit}. As the same as some old machines, let us do the subtraction without the guard digit:
$$
\begin{array}{rlcrll}
 &\;\;\;1.00 \times 2^0   &                & &\;\;\;1.00 \times 2^0 &\\
 &-1.11 \times 2^{-1} &\Longrightarrow & &-0.11  \times 2^0& \mbox{(last digit dropped)} \\
\cline{4-5}
 &                   &                & &\;\;\;0.01 \times 2^0 &= 1.00\times 2^{-2}
\end{array}
$$
The computed solution is twice, and so has relative error $100\%$.
The lack of a guard digit is a serious drawback. Almost all modern processors use the guard digit.
\vsp 

 \begin{workout}
 Assume that $x$ and $y$ are floating-point numbers in base $2$ with $y/2\leqslant x\leqslant 2y$. Show that $\fl(x-y)=x-y$ provided that the guard digit is supported and $x-y$ does not underflow.
 \end{workout}
\vsp 

Some computers perform a {\bf fused multiply-add} (FMA) operation that enables expressions like
$(x\ast y) \pm z$ (a floating-point operation followed by an addition/subtraction)
to commit just one rounding error\footnote{For example Intel Itanium,
IBM RISC System/6000 and IBM Power PC.}
\begin{equation*}
  \fl(x\ast y \pm z) = (x\ast y \pm z)(1+\delta),\quad |\delta|\leqslant u
  \vsp
\end{equation*}
for $x,y,z\in \F$. This capability enables the number of rounding
errors in many algorithms to be approximately halved. For example the
inner product $s = x^Ty$ between two vectors $x,y\in\F^n$ can be computed with just $n$ rounding
errors instead of the usual $2n-1$. See the code below.
\begin{shaded}
\vspace*{-0.5cm}
\begin{verbatim}
def InnerProd(x,y):
    s = 0
    for k in range(len(x)):
        s = x[k]*y[k] + s
    return s
\end{verbatim}
\vspace*{-0.5cm}
\end{shaded}
\noindent

\begin{workout}
Consider the Newton's method for finding the root of $f(x)=a-1/x$ for a given real number $a$.
Show that the computation of $x_{k+1}$ from $x_k$ (successive Newton's iterations)
can be expressed as two multiply-adds, thus the roundoff errors is reduced by a factor of ${1}/{4}$ if FMA operation is available.
\end{workout}
\vsp 

We close this subsection by noting that one can use the modified model
\vsp 
\begin{equation}\label{fl_op_rule2}
  \fl(x\ast y) = \frac {x\ast y}{1+\delta}, \quad |\delta|\leqslant u,\quad \ast\in\{+,-,\times,/ \}, \quad x,y\in \F
  \vsp
\end{equation}
instead of model \eqref{fl_op_rule}. 
The proof is straightforward and is left as an exercise for the reader.
For common numbers $x$ any $y$, the value $\delta$ may differ in \eqref{fl_op_rule} and \eqref{fl_op_rule2} but $|\delta|$ is bounded by the same $u$ in both models.
\vsp
\begin{workout}
Prove \eqref{fl_op_rule2}.
\end{workout}
\vsp
\section{Conditioning of problems and stability of algorithms }
Suppose we are given a mathematical problem (mathematical model) to solve.
In practice, the problem we actually solve involves inputs which are perturbed versions of the actual input data; for example, the data may be contaminated by measurement errors or roundoff errors in computer. Thus, it is important to have some knowledge about the sensitivity of the model to such perturbations in its inputs.
The fundamental question we face is as follows:
\begin{quote}
 {\em How much the solution of the perturbed problem is close to the solution of the original problem?}
\end{quote}
This concept is usually referred to as the {\bf conditioning} of the problem. A problem is called {\bf well-conditioned} if it has few sensitivity to the perturbations, i.e., the solution of the perturbed problem remains close to the solution of the original problem. Otherwise the problem is called {\bf ill-conditioned}.

We should distinguish between the conditioning of a mathematical problem and the conditioning of a computational algorithm we design to solve it. The latter is usually referred to as the {\bf stability} of the algorithm. However, sometimes the term {\em stability} is used for both mathematical model and algorithm. We discuss these concepts in the following two sections.
\vsp 
\subsection{Conditioning of a mathematical problem}
First, we try to involve you into the concept through some simple examples. Then we will investigate the theory behind our observations.
\vsp 
\begin{example}\label{ex:hilbertsys}
Many practical problems raise the need for solving a linear system of equations of the form
\begin{equation*}
  Ax = b
\end{equation*}
as a mathematical model, where $A\in \R^{n\times n}$ (i.e., an $n\times n$ matrix with real entries) and
$b\in \R^n$ (i.e., an $n$ vector with real entries) are given inputs and $x\in \R^n$ is the output (solution) we should compute. In some other applications, $A$ may be a rectangular matrix with different number of rows and columns. However, in this example we assume that $A$ is a square and nonsingular matrix. Assume further that our modeling results in a matrix $A$ of the {\em Hilbert} form
\begin{equation}\label{hilbmat}
  A = \begin{bmatrix}
        1 & \frac{1}{2} & \cdots & \frac{1}{n} \\
        \frac{1}{2} & \frac{1}{3} & \cdots & \frac{1}{n+1} \\
        \vdots & \vdots &  & \vdots \\
        \frac{1}{n} & \frac{1}{n+1} & \cdots & \frac{1}{2n-1}
      \end{bmatrix}=:H_n.
\end{equation}\\
This matrix is symmetric and positive definite and thus nonsingular. Assume that the vector $b$ is defined as the sum of all columns of $A$. Then the exact solution $x$ of $Ax=b$ is indeed $x=[1,1,\ldots,1]^T$. Let us ignore the exact solution and solve the system using a stable algorithm such as {\em Gauss elimination with pivoting}. In Python such algorithm is available via the \verb+solve+ command in the \verb+numpy+ library in the submodule \verb+linalg+. This submodule contains many other linear algebra solvers as well. We need to define a subroutine for creating a Hilbert matrix of size $n\times n$.
In the script below we assume $n=11$, form a linear Hilbert system of equations and solve it using the \verb+solve+ command from the mentioned library. Outputs are the computed solution $\wh x$ ($=$ \verb+xh+) and the relative error in the infinity norm.

\begin{shaded}
\vspace*{-3mm}
\begin{verbatim}
import numpy as np
def hilbert(n):
    A = np.zeros([n,n]);
    for i in range(n):
        for j in range (n):
            A[i,j]=1/(i+j+1)
    return A
n = 11
A = hilbert(n)
b, x = np.sum(A, axis = 0), np.ones(n)
xh = np.linalg.solve(A,b)
print('xh = ', np.round(xh,8))
e = np.linalg.norm(x-xh,np.inf)/np.linalg.norm(x,np.inf)
print('RelErr = ', np.round(e,8))
\end{verbatim}
\vspace*{-3mm}
\end{shaded}
\noindent
The output is

\begin{shaded}
\vspace*{-3mm}
\begin{verbatim}
xh =  [0.99999999 1.00000065 0.99998308 1.0001891  0.99887333 1.00396191
       0.99137162 1.0117659  0.99022415 1.00452416 0.99910608]
RelErr =  0.011765895169415952
\end{verbatim}
\vspace*{-3mm}
\end{shaded}
\noindent
We observe an unexpected result. We would typically expect a relative error close to machine precision (approximately \( u \approx 10^{-16} \) in double precision) for this small system size. However, we observe an error of the order \( 10^{-2} \) in the output. This indicates a loss of around 14 significant decimal digits in the computation; our computed result is \( 10^{14} \) times worse than what we might expect from an ideal computation. What happened? One might initially blame the Python solver \verb+solve+, but this function is based on the stable Gaussian elimination algorithm with pivoting. The true source of this serious issue lies in the original model \( Ax = b \), as we will soon see.
\end{example}
\vsp

\begin{example}[James H. Wilkinson's Example]\label{ex:polyroot}
Consider the problem of finding the roots of a polynomial of degree $n$ of the form
\begin{equation}\label{poly_n_rep}
  p(x) = x^n + a_{n-1}x^{n-1}+\cdots +a_1x +a_0.
\end{equation}
for known real coefficients $a_0,a_1,\ldots,a_{n-1}$.
If $\xi\in \C$ is a root, then $\xi$ is a function of coefficients, say $\xi:\R^n\to \C$.
In this example, we are going to test the sensitivity of $\xi$ with respect to perturbation in input values $a_0,\ldots,a_{n-1}$.
We start with
the following polynomial of degree eight\footnote{The original Wilkinson's polynomial is of degree $20$.}:
\begin{equation}\label{poly8rep}
\begin{split}
  p(x) = x^8 & -36x^7 + 546x^6 - 4536x^5+22449x^4 \\ &-67284x^3+118124x^2-109584x+40320.
\end{split}
\end{equation}
This nasty polynomial has indeed the beautiful representation $p(x)=(x-8)(x-7)\cdots(x-1)$. The exact roots are real and non-repeated numbers
$\xi = 8,7,\ldots,1$. However, in practice, a polynomial is usually represented in terms of its coefficients as shown in \eqref{poly_n_rep}, and one often needs to solve $ p(x) = 0 $ to obtain some or all of its roots.

Several approaches exist for polynomial rootfinding. A dominant algorithm involves forming the {\em companion matrix}
$$
A = \begin{bmatrix}
      -a_{n-1} & -a_{n-2} & \cdots & -a_1 & a_0 \\
      1 & 0 & \cdots & 0 & 0 \\
      0 & 1 & \cdots & 0 & 0 \\
      \vdots & \ddots & \ddots & \ddots & \vdots \\
      0 & \cdots & 0 & 1 & 0
    \end{bmatrix}\vsp
$$
where the polynomial \eqref{poly_n_rep} is its characteristic polynomial (prove this!). An iterative algorithm based on QR factorization is then applied to compute the eigenvalues of. Then an iterative algorithm based on QR factorization is applied to compute the eigenvalues of $A$ (roots of $p$). This algorithm has implemented in \verb+numpy+ library in Python. The algorithm works well for low-degree polynomials\footnote{A more efficient algorithm which does not consider the roots as a function of coefficients (but as a function of polynomial values) exists, yet not implemented in Python. Search polynomial rootfinding in {\em Chebfun}. Chebfun is an open-source software system written in MATLAB for numerical computation.}.

In the following Python code, we use the \verb+roots+ function to find the roots of \eqref{poly8rep}. To observe the sensitivity of the roots to perturbations in the coefficients, we keep all coefficients of \eqref{poly8rep} unchanged except for $a_7 = -36$, which we perturb by $0.001$. The perturbed roots are then computed using the same function.

\begin{shaded}
\vspace*{-3mm}
\begin{verbatim}
import numpy as np
coeffs = [1, -36, 546, -4536, 22449, -67284, 118124, -109584, 40320]
original_roots = np.roots(coeffs)
perturbed_coeffs = coeffs.copy()
perturbed_coeffs[1] -= 0.001
perturbed_roots = np.roots(perturbed_coeffs)
print('OriginalRoots = ', original_roots)
print('PerturbedRoots = ', perturbed_roots)
\end{verbatim}
\vspace*{-3mm}
\end{shaded}
\noindent
The output is:
\begin{shaded}
\vspace*{-3mm}
\begin{verbatim}
OriginalRoots = [8. 7. 6. 5. 4. 3. 2. 1.]
PerturbedRoots = [8.27260278  6.49985871+0.7292706j  6.49985871-0.7292706j  
                  4.57483609  4.16253083  2.99113515  2.00017793  0.9999998]
 \end{verbatim}
\vspace*{-3mm}
\end{shaded}
\noindent
As we observe, a small perturbation in one of the coefficients results in significant changes in the solution, to the extent that some roots (the second and third roots) become complex. In this example, the rootfinding algorithm performed its task correctly, and \verb+PerturbedRoots+ is indeed the vector of roots of the perturbed polynomial. This example illustrates that the roots of a polynomial, when considered {\em as functions of its coefficients}, are inherently ill-conditioned. This sensitivity to coefficient perturbations is independent of the algorithm used to calculate the roots. We will analyze this phenomenon in more detail shortly.

\end{example}
\vsp

In the above examples, we encountered several ill-conditioned problems. But how can we measure well- or ill-conditioning quantitatively? 
To have a general definition, assume that $x\in D\subset \R^m$ is the input vector, $y\in\R^m$ is the output vector, and $F$ is a map (problem model) that relates $ x $ and $ y $ via
\begin{equation}\label{F(x,y)=0}
  F(x,y) = 0.
\end{equation}
Sometimes, the output $ y $ can be explicitly represented in terms of the input $ x $. In this case, there exists a function $f:D\to\R^n$ such that
\begin{equation}\label{y=f(x)}
  y = f(x).
\end{equation}
\vsp 
\begin{definition}\label{def:conditioning}
Let \( x \in D \) be the input vector, and let \( \delta x \) be a perturbation in \( x \) such that \( x + \delta x \in D \). Similarly, let \( y \) be the output vector and \( \delta y \) be a perturbation in \( y \) such that
\[
F(x + \delta x, y + \delta y) = 0,
\]
or in explicit form,
$$
y + \delta y = f(x + \delta x).
$$
The problem given by \( F(x, y) = 0 \) or equivalently \( y = f(x) \) is said to be well-conditioned if it possesses a {\em unique solution} and the perturbations in \( y \) are bounded in a controlled manner by the perturbations in \( x \). Specifically, the condition for well-conditioning is
\vsp
$$
\frac{\|\delta y\|}{\|y\|} \leq C \frac{\|\delta x\|}{\|x\|}
$$
for \( x \neq 0 \) and \( y \neq 0 \), where \( C \) is a relatively small constant.
In special cases when \( x = 0 \) and \( y \neq 0 \) the above bound is replaced by
\vsp 
   \[
   \frac{\|\delta y\|}{\|y\|} \leq C \|\delta x\|.
   \]
If \( x \neq 0 \) and \( y = 0 \) we may write
   \[
   \|\delta y\| \leq C \frac{\|\delta x\|}{\|x\|},
   \]
and finally if \( x = 0 \) and \( y = 0 \) then we must replace it by 
   \[
   \|\delta y\| \leq C \|\delta x\|.
   \]
In each of these cases, the constant \( C \) reflects how sensitively the solution \( y \) responds to changes or perturbations in the input \( x \). A small \( C \) indicates that the problem is well-conditioned, meaning small changes in the input result in small changes in the output. Conversely, if \( C \) is large, the problem is ill-conditioned, implying that small changes in the input can cause large changes in the output, which can be problematic for numerical computations.
\end{definition}
\vsp 

For many problems, the constant $C$ in the above definition can be estimated. This estimate of $C$ is known as the condition number of the problem. The condition number provides a measure of how sensitive the solution of a problem is to changes or errors in the input.
According to Definition \eqref{def:conditioning},
the condition number of a problem can be obtained as the following ratio:\\
 \[
 \frac{\mbox{amount of perturbation in the output (solution) of the problem}}{\mbox{amount of perturbation in the input of the problem}}
 \]
\ \\
Depending on the criterion used to measure the amount of perturbation in the input and in the output of a given problem, {\em relative} or {\em absolute} condition numbers are defined.

Consider the simplest case $m=n=1$ with the explicit representation \eqref{y=f(x)} for the problem. Assume that $f$ is twice continuously differentiable. The Taylor expansion then gives
\vsp
$$
y + \delta y = f(x+\delta x) = f(x) + \delta x f'(x) + \mathcal O(|\delta x|^2).
$$
Ignoring the error term for small input perturbation $\delta x$, we can write
\vsp 
$$
\frac{\delta y}{y} \approx \frac{xf'(x)}{f(x)}\cdot \frac{\delta x}{x}.
\vsp
$$
provided that $x\neq 0$ and $y\neq 0$. Consequently, we can define the relative {\em condition number of} $f$ {\em at point} $x$ by
\begin{equation}\label{cond:def1}
  (\cond\, f)(x) := \frac{|x||f'(x)|}{|f(x)|}
\end{equation}
to have
$$
\frac{|\delta y|}{|y|}\approx (\cond f)(x) \frac{|\delta x|}{|x|}.
$$
\\
The condition number shows, approximately, how much the perturbation in the input data amplifies in the solution.
A large condition number indicates ill-conditioning, while a small condition number signifies well-conditioning. The threshold for what constitutes ``large'' or ``small'' depends on the specific problem and the desired accuracy. If \( x = 0 \) but \( y \neq 0 \), then \(\delta x\) should be measured absolutely and \(\delta y\) relatively. In this case we define
\vsp 
\begin{equation*}
  (\cond\, f)(x) := \frac{|f'(x)|}{|f(x)|}, \quad \frac{|\delta y|}{|y|}\approx (\cond f)(x) {|\delta x|}.
\end{equation*}
On the other hand, if $x\neq0$ and $y=0$, then
\begin{equation*}
  (\cond\, f)(x) := {|x||f'(x)|}, \quad {|\delta y|}\approx (\cond f)(x) \frac{|\delta x|}{|x|}.
\end{equation*}
Finally, if $x=y=0$ the condition number of $f$ can be defined as 
\vsp
\begin{equation*}
  (\cond\, f)(x) := {|f'(x)|}, \quad {|\delta y|}\approx (\cond f)(x) {|\delta x|}.
\end{equation*}
\vsp 
\begin{example}
Consider the solution $y$ of the quadratic equation $y^2-2ay+1=0$ for input value $a>1$.
In the implicit form we may write $F(a,y)=y^2-2ay+1=0$, but solving $y$ in terms of $x$ gives the explicit form
$$
y_{\pm} = a \pm \sqrt{a^2-1}=:f_{\pm}(a).
$$
This problem indeed has two solutions. We can treat \( y_{+} \) and \( y_{-} \) either individually or together. Since $f'_{\pm}(a)=1\pm \frac{a}{\sqrt{a^2-1}}$, from \eqref{cond:def1} we have
$$
(\cond\, f_{\pm})(a) := \frac{|f'_{\pm}(a)||a|}{|f_{\pm}(a)|} = \frac{|a|}{\sqrt{a^2-1}}, \quad a>1.
$$
This shows that for values of $a$ far from $1$, the problem is well-conditioned, while for values of $a$ close to $1$ (i.e. when the roots tend to become of multiplicity $2$) the problem becomes ill-conditioned.

The ill-conditioning can be bypassed by using a simple change of variable to obtain an equivalent but well-conditioned problem for values of \(a\) near $1$, and even for \(a = 1\).
If we use $b=a+\sqrt{a^2-1}$ then the quadratic equation is reformulated as
\vsp 
$$
y^2 -\frac{1+b^2}{b}y +1 = 0.
$$
In this case we have
$$
y_{+}=f_{+}(b)=\frac{1}{b}, \quad y_{-}=f_{-}(b) = b.
$$
It is left to you to show that both $f_{+}$ and $f_{-}$ are well-conditioned functions for values of $b$ (or $a$) close to $1$.
\end{example}
\vsp
\begin{example}

In this example, we determine the conditioning of a problem that cannot be easily represented in an explicit form. Consider the nonlinear equation
$$
x^n-ae^{-x}=0, \quad a>0, \quad n\geqslant 1.
$$
Here, we assume that \( n \) is a fixed positive integer, \( a \) is the input data, and \( x \) is the output (noting that this is contrary to our usual notation where \( x \) typically represents the input). This equation has exactly one positive root, denoted by \(\xi\) (prove this!). We aim to measure the sensitivity of this root with respect to a perturbation in \(a\). Therefore, we consider \(\xi\) as a function of \(a\), say \(\xi(a)\). By estimating the condition number \((\cond\, \xi)(a)\), we will demonstrate that \(\xi(a)\) is a well-conditioned function of \(a\).
Since an explicit form for \(\xi\) is not available, implicit differentiation can be used to compute \(\xi'(a)\).
We have
$$
[\xi(a)]^n-ae^{-\xi(a)}=0.
$$
Implicit differentiation with respect to $a$ yields
$$
n \xi'(a) [\xi(a)]^{n-1} - e^{-\xi(a)} + a\xi'(a) e^{-\xi(a)} = 0,
$$
or
$$
\xi'(a) = \frac{e^{-\xi(a)}}{n[\xi(a)]^{n-1}+ae^{-\xi(a)}}.
$$
\vsp 
Using this in the definition of condition number gives
\vsp 
$$
(\cond\, \xi)(a)=\frac{\xi'(a)a}{\xi(a)} = \frac{ae^{-\xi}}{n\xi^{n}+a\xi e^{-\xi}}=\frac{ae^{-\xi}}{nae^{-\xi}+a\xi e^{-\xi}}=\frac{1}{n+\xi}\leqslant \frac{1}{n}.
$$
Since all involved quantities are positive, absolute values were not needed in the definition of the condition number. The derived bound for the condition number indicates that the root \(\xi\) as a function of \(a\) is well-conditioned.
\end{example}
\vsp

For the case of arbitrary $m$ and $n$, assume that
$$
x = (x_1,\ldots,x_m)^T\in \R^m, \quad y=(y_1,\ldots,y_n)^T\in \R^n,
$$
and $f:\R^m\to \R^n$ maps data $x$ into solution $y$. In an element-wise  form we can write
$$
y_k =  f_k(x_1,\ldots,x_m),\quad k=1,2,\ldots,n.
$$
Assume further that each function $f_k:\R^m\to\R$ has partial derivatives with respect to all $m$ variables at point $x$.

One way to measure the sensitivity of solution $y$ with respect to small changes in $x$ is to subject only one variable, $x_j$, to a perturbation and observe the resulting change in
just one component $y_k$. Then we can apply the univariate definition \eqref{cond:def1} and obtain
\begin{equation}\label{cond:defkj}
  \kappa_{kj}(x):= \frac{\left|\frac{\partial f_k}{\partial x_j}(x)\right||x_j|}{|f_k(x)|}.
\end{equation}
If a component of $x$, or of $y$, vanishes, one should modify \eqref{cond:defkj} as discussed earlier.
This gives us a whole matrix
$$
K(x) = [\kappa_{kj}(x)] \in \R^{n\times m}
$$
of condition numbers. If a single condition number is sought, we can use a matrix norm to obtain
\begin{equation}\label{cond:def2}
  (\cond\, f)(x):=\|K(x)\|.
\end{equation}
The condition number can be estimated in another way, often simpler but sometimes misleading.
Let the relative perturbation in input $x\in\R^m$ and output $y\in\R^m$ be measured by
$$
\frac{\|\delta x\|_{\infty}}{\|x\|_{\infty}} \quad \mbox{and} \quad \frac{\|\delta y\|_{\infty}}{\|y\|_{\infty}},
$$
respectively, where $\delta x=(\delta x_1,\ldots, \delta x_m)^T$ and $\delta y=(\delta y_1,\ldots, \delta y_n)^T$. Now, in analogy to the scalar case and using the linear multivariate Taylor expansion, we have
\vsp 
$$
y_k + \delta y_k = f_k(x+\delta x) = f_k(x)+ \sum_{j=1}^{m}\frac{\partial f_k}{\partial x_j}(x) \delta x_j + \mathcal{O}(\|\delta x\|_\infty^2).
$$
Ignoring the error term for a small input perturbation $\delta x$, we have, at least approximately,
\vsp 
$$
|\delta y_k|\leqslant \sum_{j=1}^{m}\left|\frac{\partial f_k}{\partial x_j}(x)\right| |\delta x_j|
\leqslant \|\delta x\|_\infty \max_{k}\sum_{j=1}^{m}\left|\frac{\partial f_k}{\partial x_j}(x)\right|
$$
for all $k=1,2,\ldots,n$. This simply gives\footnote{Remind that the norm infinity of a $n\times m$ matrix $A$ is defined as
$\ds \|A\|_\infty = \max_{1\leqslant k\leqslant n}\sum_{j=1}^{m}|a_{kj}|$.}
\vsp
\begin{equation}\label{delylessdelx}
\|\delta y\|_\infty \leqslant \|\delta x\|_\infty \|J_f(x)\|_\infty,
\end{equation}
where $J_f(x)$ is the {\em Jacobian matrix} defined by
\vsp 
$$
J_f(x) = \begin{bmatrix}
           \frac{\partial f_1}{\partial x_1} & \frac{\partial f_1}{\partial x_2} & \cdots & \frac{\partial f_1}{\partial x_m} \\
           \vdots & \vdots &  & \vdots \\
           \frac{\partial f_n}{\partial x_1} & \frac{\partial f_n}{\partial x_2} & \cdots & \frac{\partial f_n}{\partial x_m}
         \end{bmatrix}\in \R^{n\times m}.
$$
For $x\neq 0$ and $y\neq 0$ from \eqref{delylessdelx} we have
\vsp 
\begin{equation*}
  \frac{\|\delta y\|_\infty}{\|y\|_\infty}\leqslant \frac{\|J_f(x)\|_\infty\|x\|_\infty}{\|f(x)\|_\infty}\cdot \frac{\|\delta x\|_\infty}{\|x\|_\infty}
\end{equation*}
which suggests us to define
\begin{equation}\label{cond:def3}
  (\cond\,f)(x):= \frac{\|J_f(x)\|_\infty\|x\|_\infty}{\|f(x)\|_\infty}.
\end{equation}
In situations where either $x=0$ or $y=0$, a modification in definition is carried out similar to the univariate case.
For $m=n=1$, this definition of condition number reduces to \eqref{cond:def1}.
However, for $m,n>1$ the condition number in \eqref{cond:def3} may mislead us from the actual conditioning of the problem.
The reason is that the norms sometimes tend to destroy the details. For example, if $x$ has components
of vastly different magnitudes, then $\|x\|_\infty$ is simply equal to the largest of these
components, and all the others are ignored. See the workout below for an example.
\vsp 
\begin{workout}
Let $x=[x_1,x_2]^T$ and
$$
f(x) = \begin{bmatrix}
         \frac{1}{x_1}+\frac{1}{x_2} \\
         \frac{1}{x_1}-\frac{1}{x_2}
       \end{bmatrix}.
$$
Compute the condition number of $f$ using both formulas \eqref{cond:def2} and \eqref{cond:def3}. Show that the formula 
\eqref{cond:def2} exhibits the potential for ill-conditioning for certain values of $x$ while the formula \eqref{cond:def3} indicates that  $f$ is a well-conditioned function for all values of $x$. 	
\end{workout}
\vsp 
\begin{example}
Coming back to Example \eqref{ex:polyroot}, here we analyze the ill-conditioning of the roots of a polynomial as functions of its coefficients. Let $\xi$ be a fixed root of polynomial
\begin{equation*}
  p(x) = x^n + a_{n-1}x^{n-1}+\cdots +a_1x +a_0
\end{equation*}
where $a_0\neq 0$, i.e., $\xi\neq 0$. Besides, assume that $\xi$ is a simple root, i.e.,
$p'(\xi)\neq 0$. The root $\xi$ is a complex function of coefficients $a=(a_{0},\ldots,a_{n-1})$;
$$
\xi:\R^n \to \C,\quad \xi = \xi(a_{0},\ldots,a_{n-1}).
$$
Recalling \eqref{cond:defkj}, we have
\begin{equation}\label{cond:kappa_k}
\kappa_k(a) = \frac{\left|\frac{\partial \xi}{\partial a_k}(a)\right||a_k|}{|\xi|}.
\end{equation}
The vector of condition numbers then is $K(a)=(\kappa_0,\ldots,\kappa_{n-1})$.
To have an explicit representation for $\kappa_k$ we must compute the partial derivatives of $\xi$ with respect to coefficients $a_k$. Since $\xi$ is a root, we have
$$
\xi^n + a_{n-1}\xi^{n-1}+\cdots+a_k\xi^k + \cdots + a_1\xi+a_0 = 0.
$$
An implicit differentiation with respect to $a_k$ yields
\vsp 
$$
n\xi^{n-1}\frac{\partial \xi}{\partial a_k} + (n-1)\xi^{n-2}\frac{\partial \xi}{\partial a_k}+\cdots+ \xi^k + a_k\xi^{k-1}\frac{\partial \xi}{\partial a_k} +\cdots+ a_1\frac{\partial \xi}{\partial a_k}+ 0 \equiv 0.
$$
This equation is equivalent to
$$
p'(\xi) \frac{\partial \xi}{\partial a_k} +\xi^{k} = 0.
$$
Since $p'(\xi)\neq0$, we can write
$$
\frac{\partial \xi}{\partial a_k} = -\frac{\xi^k}{p'(\xi)}.
$$
\vsp
Inserting into \eqref{cond:kappa_k} and using the $\|\cdot\|_1$ for condition number, we obtain
\begin{equation}\label{cond:pkappa}
(\cond\, \xi)(a) = \|K(x)\|_1 = \frac{1}{|\xi p'(\xi)|}\sum_{k=0}^{n-1}|a_k||\xi|^{k}=\frac{\ds\sum_{k=0}^{n-1}|a_k||\xi|^{k}}{\left|\ds\sum_{k=0}^{n-1}ka_k\xi^{k}\right|}.
\end{equation}
The formula \eqref{cond:pkappa} shows that the problem has a potential to become ill-conditioned as (depending on signs of $a_k$ and $\xi$) the denominator could be much smaller than the numerator. As an example, the condition numbers for roots of polynomial
\vsp 
\begin{align*}
  p(x) = &\, x^8 -36x^7 + 546x^6 - 4536x^5+22449x^4-67284x^3+118124x^2-109584x+40320\\
       =& \,(x-8)(x-7)(x-6)(x-5)(x-4)(x-3)(x-2)(x-1),
\end{align*}
considered in Example \ref{ex:polyroot}, are computed using formula \eqref{cond:pkappa} and given in Table \ref{tb-cond1-poly}.

\begin{center}
\captionof{table}{Condition numbers of the roots of polynomial \eqref{poly8rep}}\label{tb-cond1-poly}
\begin{tabular}{c|cccc}
\hline
  $\xi_j$ & $8$  & $7$   &  $6$ & $5$    \\
  \hline
$(\cond\, \xi_j)(a)$ & $0.14\times 10^{9}$ & $0.47\times 10^{9}$  & $0.44\times 10^{9}$ & $0.18\times 10^{9}$ \\
  \hline
    \hline
      $\xi_j$ & $4$  & $3$   &  $2$ & $1$    \\
      \hline
$(\cond\, \xi_j)(a)$ & $0.35\times 10^{8}$ &$0.25\times 10^{7}$   & $0.48\times 10^{5}$& $0.72\times 10^{2}$ \\
  \hline
\end{tabular}
\end{center}
Once more, look at the output \verb+PerturbedRoots+ in the Python code provided in Example \ref{ex:polyroot} to observe that the larger the condition number of a root, the more significant the change in the root becomes. 
The worst condition numbers are associated with roots $\xi=6, 7$ which result in a pair of complex perturbed roots.
\end{example}
\vsp 
\begin{example}[Conditioning of a Linear System of Equations]\label{ex:condA}
Here we analyze the disaster observed in Example \ref{ex:hilbertsys} using the concept of condition numbers. Consider again the linear system
\begin{equation}\label{sys_Axb}
Ax = b,
\end{equation}
for a given nonsingular $n\times n$ matrix $A$ and a nonzero vector $b\in \R^n$. Here, $x\in \R^n$ is served as the solution (output) and we aim to investigate the conditioning of the problem when the input data $A$ and $b$ are subjected to small perturbations. For simplicity, we only perturb the vector $b$ and keep the matrix $A$ unchanged. So, we consider $x$ as a function of $b$ only;
\vsp
$$
x = A^{-1}b := f(b).
$$
Since $J_f(b)=A^{-1}$, from \eqref{cond:def3} we have
\begin{equation*}
  (\cond\, f)(b) = \frac{\|b\|\|A\|}{\|A^{-1}b\|}=\frac{\|Ax\|\|A^{-1}\|}{\|x\|},
\end{equation*}
and since there is a one-to-one correspondence between $x$ and $b$, we find for the
worst condition number
$$
\max_{b\neq 0}\, (\cond\, f)(b) = \max_{x\neq 0}\frac{\|Ax\|}{\|x\|}\cdot \|A^{-1}\| = \|A\|\cdot \|A^{-1}\|
$$
using the definition of a natural norm of $A$. The number on the far right no longer depends on $b$ and is called the {\bf condition number} of the matrix
$A$. We denote it by
\begin{equation}\label{condA}
  \cond( A) := \|A\|\cdot \|A^{-1}\|.
  \vsp
\end{equation}
To get a deeper insight into $\cond(A)$, consider the perturbed system
$$
A(x+\delta x) = b+\delta b
$$
which together with \eqref{sys_Axb} gives the explicit representation $\delta x = A^{-1}\delta b$ for the output perturbation $\delta x$. Taking norm from both sides of this relation gives
$$
\|\delta x\| \leqslant \|A^{-1}\| \|\delta b\|.
$$
On the other side, from original system \eqref{sys_Axb} we have $\|b\|\leqslant \|A\|\|x\|$, or
$$
\frac{1}{\|x\|}\leqslant \frac{\|A\|}{\|b\|}.
$$
Multiplying both sides of recent equations, we get
\vsp
\begin{equation}\label{cond:inlinsys}
\frac{\|\delta x\|}{\|x\|}\leqslant \|A\|\cdot \|A^{-1}\|\frac{\|\delta b\|}{\|b\|}= \cond(A)\frac{\|\delta b\|}{\|b\|}.
\end{equation}
The error bound \eqref{cond:inlinsys} simply describes what happened in Example \ref{ex:hilbertsys} for the Hilbert system.
In fact,
rounding errors imply a small perturbation in the input data of machine epsilon order, i.e.,
$$
\frac{\|\delta b\|_\infty}{\|b\|_\infty} = u \approx 10^{-16}
$$
in the double precision floating-point format. The condition number of a Hilbert matrix of size $11\times 11$ is
approximately $10^{+15}$. The error bound \eqref{cond:inlinsys} then gives
\vsp 
$$
\frac{\|x-\wh x\|_\infty}{\|x\|_\infty}\lesssim 10^{+15}\cdot 10^{-16} = 0.1
$$
which indicates that approximately 15 decimal digits will be lost when solving such a Hilbert system with any standard algorithm in double precision format. Our observation in Example \ref{ex:hilbertsys}
confirms this conclusion.
\end{example}
\vsp

Although in Example \eqref{ex:condA} we have considered only perturbations in the right-hand vector $b$, it
turns out that the error bounds still depend on the condition number in \eqref{condA} when we also account for perturbations in the matrix $A$. See Workout \ref{wo1:condA} below.
\vsp

\begin{workout}\label{wo1:condA}
Consider the linear system of equations $Ax=b$ for nonsingular matrix $A\in \R^{n\times n}$. Let $\delta x$ be the perturbation raised in $x$ caused by changing $A$ to $A+\delta A$ but keeping $b$ unchanged. Prove that
\begin{equation*}
\frac{\|\delta x\|}{\|x\|}\leqslant \frac{\cond(A)\frac{\|\delta A\|}{\|A\|}}{1-\cond(A)\frac{\|\delta A\|}{\|A\|}},
\end{equation*}\\
provided that $\delta A$ is so small such that $\|\delta A\|\|A^{-1}\|\leqslant 1$. Moreover, prove that under perturbation in both $A$ and $b$ we have
\vsp 
\begin{equation*}
\frac{\|\delta x\|}{\|x\|}\leqslant \frac{\cond(A)}{1-\cond(A)\frac{\|\delta A\|}{\|A\|}}\left( \frac{\|\delta b\|}{\|b\|} + \frac{\|\delta A\|}{\|A\|} \right).
\end{equation*}
\end{workout}
\vsp 

It is interpreted from Example \ref{ex:condA} and Workout \ref{wo1:condA} that the condition number of a matrix plays a crucial role in numerical matrix computations. 

As an example, the Hilbert matrix \eqref{hilbmat} is an ill-conditioned matrix. Table \ref{tb:condhilb} shows the condition number of this matrix for different values of \( n \) in both the infinity norm and the 2-norm. We observe a rapid increase in the condition numbers as \( n \) grows. A Hilbert matrix with small size $13\times 13$, which has a condition number of order $10^{17}$, destroys all $16$ decimal significant digits in the double precision floating-point format.

\begin{table}[th!]
  \centering
  \caption{Condition numbers of the Hilbert matrices of different sizes}\label{tb:condhilb}
  \begin{tabular}{c|cccccc}
    \hline
    $n$ &          $3$ &  $5$ & $7$& $9$ & $11$ & $13$ \\ \hline \\
    $\cond_2(H_n)$    & $5.24\times 10^{2}$  & $4.77\times 10^{5}$ & $4.75\times 10^{8}$ & $4.93\times 10^{11}$ & $5.22\times 10^{14}$ & $4.79\times 10^{17}$\\\\
    $\cond_\infty(H_n)$& $7.48\times 10^{2}$ & $9.44\times 10^{5}$ & $9.85\times 10^{8}$ & $1.10\times 10^{12}$ & $1.23\times 10^{15}$ & $8.53\times 10^{17}$\\
    \hline
  \end{tabular}
\end{table}

Another well-known ill-conditioned matrix is the {\em Vandermonde matrix}. This matrix is of the form
\begin{equation*}
  V_n = \begin{bmatrix}
          1 & 1 & \cdots & 1 \\
          t_1 & t_2 & \cdots & t_n \\
          \vdots & \vdots &  & \vdots \\
          t_1^{n-1} & t_2^{n-1} & \cdots & t_n^{n-1}
        \end{bmatrix}\in \R^{n\times n}
\end{equation*}\\
where $t_1,t_2,\ldots,t_n$ are some distinct real numbers. For example, for equality spaced numbers 
\vsp
$$
t_k = -1+\frac{2(k-1)}{n-1}, \quad k=1,2,\ldots, n,
$$
in interval $[-1,1]$, the condition numbers of $V_n$ are given in Table \ref{tb:condvander} for some values of $n$.
Although they do not grow as quickly as those for the Hilbert matrix, they still increase exponentially fast. Worse than exponential growth is observed if one takes harmonic numbers
\vsp
$$
t_k = \frac{1}{k}, \quad k=1,2,\ldots n.
$$
In this case, we can show that the condition number of \( V_n \) grows as \( n^{n+1} \), which is significantly worse than the condition number of the Hilbert matrix \( H_n \).

\begin{table}[th!]
  \centering
  \caption{Condition numbers of the Vandermonde matrices of different sizes}\label{tb:condvander}
  \begin{tabular}{c|cccccc}
    \hline
    $n$ &          $10$ &  $15$ & $20$& $25$ & $30$ & $35$ \\ \hline \\
    $\cond_2(V_n)$    & $2.63\times 10^{3}$  & $1.10\times 10^{6}$ & $2.72\times 10^{8}$ & $7.05\times 10^{10}$ & $1.84\times 10^{13}$ & $4.61\times 10^{15}$\\\\
    $\cond_\infty(V_n)$& $2.06\times 10^{4}$ & $5.58\times 10^{6}$ & $1.75\times 10^{9}$ & $4.91\times 10^{11}$ & $1.46\times 10^{14}$ & $4.16\times 10^{16}$\\
    \hline
  \end{tabular}
\end{table}

In scientific computing, it is important to avoid mathematical and computational models that result in ill-conditioned matrices whenever possible. There are often alternative models and simulations that lead to well-conditioned matrices. 
For example, when solving polynomial interpolation using the monomial basis \(\{1, x, \ldots, x^n\}\), the algorithm involves the Vandermonde matrix, which can be ill-conditioned. This can be circumvented using other algorithms such as Newton's method\footnote{Newton's method may introduce other instability issues.} or the {\em  barycentric Lagrange interpolation} method.
Similarly, solving the best polynomial approximation problem in the 2-norm leads to a Hilbert system if the monomial basis \(\{1, x, \ldots, x^n\}\) is used. This issue can be avoided by employing an orthogonal basis instead.

When there is no flexibility to choose alternative models and methods, matrix {\em preconditioners} can be helpful. If \( A \) is an ill-conditioned matrix, a preconditioner \( P \) is a nonsingular matrix such that \( PA \) is well-conditioned. Then we can solve 
$$
PAx = Pb
$$
instead of original system $Ax = b$. However, constructing a preconditioner \( P \) for a matrix \( A \) is not straightforward in many circumstances. There is no universal approach for building preconditioners that work effectively for all matrices.
\vsp

\begin{remark}
The condition number of a square, nonsingular matrix \( A \) is defined by \eqref{condA}. In a forthcoming lecture, we will explore an extension of this concept to non-square and even rank-deficient matrices using singular value decomposition (SVD).
\end{remark}
\vsp
\begin{remark}
In Python, an algorithm for estimating the condition numbers of a matrix is implemented in the \verb+numpy.linalg+ library. It works for $1$, $2$, infinity (\verb+np.inf+) and Frobenius (\verb+'fro'+) norms. For instance, to get the condition number in norm infinity we write
\begin{verbatim}
       np.linalg.cond(A,np.inf)
\end{verbatim}
The default command 
\verb+np.linalg.cond(A)+ gives the condition number in the Frobenius norm.
\end{remark}
\vsp

\subsection{Stability of an algorithm}

There might be different numerical algorithms for solving a given mathematical problem. Some algorithms are more sensitive than others to small errors in the input, meaning that small perturbations might cause large perturbations in the output of the algorithm.

The precise definition of the concept of stability varies depending on the subject area. While we do not intend to provide a rigorous definition, we aim to offer a general insight into the concept of stability through a definition and some concrete examples.

Mathematical models are typically approximated via a {\em numerical method} by converting to a new {\bf discretized} version with possibly new approximated inputs. The solution of the discretized problem is expected to well approximate the solution of the original problem. Analogous to the primary problem \eqref{F(x,y)=0} or \eqref{y=f(x)}, the discretized problem can be expressed as
\vsp
\begin{equation}\label{FN(xN,yN)=0}
F_N(x_N,y_N) = 0, \quad \text{or} \quad y_N = f_N(x_N)
\end{equation}
where \( F_N \) represents the discretized model, \( x_N \) denotes approximate inputs, and \( y_N \) represents the solution of the discretized problem. The subscript \( N \) is a discretization parameter. This is an intermediate (and crucial) step of {\em problem-solving} in scientific computing. All procedures for solving continuous problems such as ordinary and partial differential equations (ODEs and PDEs) necessarily involve this step, and a significant amount of research in numerical analysis and scientific computing is dedicated to designing and developing efficient discretization techniques for various types of mathematical problems.

The concept of {\bf stability} can be extended here in a manner similar to that described for the conditioning of a mathematical problem in the preceding section, by replacing \( f \), \( x \), and \( y \) with \( f_N \), \( x_N \), and \( y_N \), respectively. However, it is important to note that, rather than stability, the discretized problem \eqref{FN(xN,yN)=0} is subject to another issue that requires the discretized model \( F_N \) to accurately approximate the exact model \( F \) as the discretization becomes finer. This property is known as {\bf consistency}. In many cases, consistency, along with some form of stability, implies the {\bf convergence} of the approximate solution \( y_N \) to the exact solution \( y \).

In a forthcoming lecture, we will learn these concepts (consistency and stability) in detail for the numerical solution of ODEs. However, in this section (and the next section), we will assume that we are given a discretized problem to program into the computer, and our problem is only subject to roundoff errors in its inputs. Our goal is to investigate the sensitivity of the algorithm with respect to input perturbations. Roughly speaking, we refer to an algorithm as {\bf  numerically stable} if rounding errors occurring during the course of the algorithm are not amplified rapidly; in other words, if small perturbations in the inputs of the algorithm result in only small perturbations in its output.

In the stability analysis of an algorithm, different strategies can be employed, including {\bf forward analysis} and {\bf backward analysis}.
To simplify notation, we drop the subscript \( N \) from the discretization problem \eqref{FN(xN,yN)=0} and assume that
$$
F(x,y) = 0
$$
is our numerical algorithm with exact input \( x \) and exact output \( y \). The algorithm \( F \) accepts an input data \( x \) contaminated with roundoff error \( \delta \), i.e., \( x + \delta \), and produces a numerical solution \( \wh y \), which we hope to be close to the exact solution \( y \). Thus, we may write
\begin{equation}\label{stab:F}
F(x + \delta, \wh y) = 0.
\end{equation}
In a forward analysis, we assume that the relative error in \( \delta \) is proportional to the machine's precision, say \( \|\delta\|/\|x\| \leqslant C_{in} u \) where \( C_{in} \) is a small constant, and provide a bound on the error in the solution. That is, we look for a constant \( C_{out} \) such that
\vsp
\[
\frac{\|y - \wh y\|}{\|y\|} \leqslant C_{out} u.
\]
If, depending on the problem and the expected accuracy, \( C_{out} \) is a rather small constant, we say that the algorithm is {\em forward stable}.

In backward analysis, given a certain computed solution \( \wh y \), we look for a perturbation \( \delta \) on the data such that \eqref{stab:F} is satisfied. We say the algorithm is {\em backward stable} if the bound on perturbation \( \delta \) is small, i.e.,
\[
\frac{\|\delta\|}{\|x\|} \leqslant C_{bkw} u
\vsp
\]
where \( C_{bkw} \) is a small constant. In other words, an algorithm is backward stable if {\em the computed solution \( \wh y \) is the exact solution of a nearby problem.}
\vsp 

\begin{center}
\begin{tikzpicture}
    \draw[color=black,thick] (0,0)--(4,0)  node[right] {$~y$};
    \draw[color=black,thick] (1,-1)--(5,-1.7)  node[right] {$~\wh y$};
    \filldraw  (0,0) circle (2pt);
    \filldraw  (4,0) circle (2pt);
    \filldraw  (1,-1) circle (2pt);
    \filldraw  (5,-1.7) circle (2pt);
    \node[color=black] at (-0.3,0) {$x$};
    \node[color=black] at (0.5,-1.3) {$x+\delta$};
    \node[color=black] at (2,0.25) {$F$};
    \node[color=black] at (3,-1.1) {$F$};
    \draw [<->] (0.05,-0.05)--(0.95,-0.95) node[right] {};
    \draw [<->] (4.05,-0.05)--(4.95,-1.65) node[right] {};
    \node[color=black] at (5.7,-0.7) {\footnotesize{Forward error}};
    \node[color=black] at (-0.9,-0.5) {\footnotesize{Backward error}};
\end{tikzpicture}
\end{center}

Forward and backward analyses are two different instances of the so
called {\em a priori analysis}, and can be applied to investigate not only
the stability of an algorithm but also the convergence of the solution of a discretized problem (numerical method) to the solution of a mathematical model. In this case
it is referred to as {\em a priori error analysis}, which can again be performed
using either a forward or a backward technique.
In {\em a posteriori analysis}, we aim to provide an estimate on the error of the solution in terms of the reminder.
More precisely, we bound the error $y-\wh y$ as a function of the residual
$$
r = F(x,\wh y).
$$

\begin{example}
Let's revisit the linear system of equations \( Ax = b \) for a nonsingular matrix \( A \in \R^{n \times n} \) with a fixed positive integer \( n \). Here, \( A \) and \( b \) are inputs, and \( x \) is the solution. Suppose this system is solved using a numerical algorithm for linear systems of equations, such as the Gauss elimination algorithm with pivoting. In almost all cases, \( A \) and \( b \) are subject to roundoff errors of magnitudes, at least at the level of the machine's precision \( u \). Let the computed solution be denoted by \( \wh x \).

In the forward analysis for the Gauss elimination algorithm, we estimate the error \( x - \wh x \) in terms of bounds on relative perturbations in \( A \) and \( b \), here \( u \). We seek a constant \( C_{out} \) (which naturally depends on \( n \)) such that
\vsp
$$
\frac{\|x-\wh x\|}{\|x\|}\leqslant C_{out} u.
$$

In the backward analysis, on the other hand, we estimate the perturbations \( E \) and \( e \) that need to be introduced to \( A \) and \( b \), respectively, in order to obtain
$$
(A+E)\wh x  = b+e
$$
for the computed solution \( \wh x \). The bounds on \( E \) and \( e \) can be expressed in terms of \( u \):
\vsp
$$
\frac{\|E\|}{\|A\|}\leqslant C_{bkw} u, \quad \frac{\|e\|}{\|b\|}\leqslant C_{bkw} u.
$$
The constants \( C_{bkw} \) depend on \( n \), of course. Experimental evidence suggests that for most matrices appeared in practical problems, the constant \( C_{out} \) in the forward analysis and the constant \( C_{bkw} \) in the backward analysis are both of the order \( n \), the size of the system. This would constitute an acceptable bound, making the Gauss algorithm known as a stable algorithm for a majority of practical problems.

Finally, in {\em a posteriori} error analysis, we seek an estimate for the error \( x - \wh x \) as a function of the residual \( r = b - A\wh x \). Since \( x - \wh x = A^{-1}r \), we can write \( \|x - \wh x\| \leqslant \|A^{-1}\| \|r\| \). On the other hand, from \( Ax = b \), we have \( 1/\|x\| \leqslant \|A\|/\|b\| \). These two bounds together yield
\vsp
$$
\frac{\|x-\wh x\|}{\|x\|}\leqslant \frac{\cond(A)}{\|b\|}\|r\|,
$$
which provides a bound on the forward error in terms of the residual \( r \). This error bound shows that a small residual does not necessarily correspond to a small error in the solution, especially if the matrix \( A \) is ill-conditioned.

\end{example}
\vsp

At the end of this section let us address a key lesson concerning the question of
{\em when can we expect a numerical solution to be accurate?} The answer is remarked below.
\vsp
\begin{remark}
In general, we can have high confidence in the accuracy of a numerical solution if the problem is well-conditioned {\em and} a numerically stable algorithm has been employed to solve it. However, if the problem is ill-conditioned or if a numerically unstable algorithm has been used, the computed solution may be inaccurate.
\end{remark} 
\vsp

\section{Roundoff error analysis for some simple algorithms}\label{sect-analysis_algorithms}

According to \eqref{fl_op_rule}, a single floating-point operation introduces an amplification factor $(1+\delta)$ (or $(1+\delta)^{-1}$ in \eqref{fl_op_rule2}) in the computed result. For an arithmetic with multiple floating-point operations (for example summing or multiplying $n$ numbers), the expression of the result may contain a sort of products of $(1+\delta_k)^{\pm1}$ terms. The following Lemma proposes an elegant way  to simplify such expressions \cite{Higham:2002}.
\begin{lemma}\label{lem:gamma_n}
If $|\delta_k|\leqslant u$ and $\rho_k=\pm1$ for $k=1,\ldots,n$, and $nu\leqslant 1$ then
$$
\prod_{k=1}^n (1+\delta_k)^{\rho_k} = (1+\theta_n), \quad |\theta_n|\leqslant \frac{nu}{1-nu}=:\gamma_n.
$$
\end{lemma}
\noindent
\proof
We use an induction on $n$. Assume that the result is true for $n$. We have for case $\rho_{n+1}=1$,
$$
\prod_{k=1}^{n+1} (1+\delta_k)^{\rho_k}= (1+\theta_n)(1+\delta_{n+1}) =: (1+\theta_{n+1})
$$
where $\theta_{n+1}=\theta_n+\delta_{n+1}+\theta_n\delta_{n+1}$. Thus, we can write 
\vsp
\begin{align*}
  |\theta_{n+1}|\leqslant &\, |\theta_n|+|\delta_{n+1}|+|\theta_n\delta_{n+1}|\\
                \leqslant &\, \gamma_n + u + u\gamma_n  = \frac{nu}{1-nu} + u + \frac{nu^2}{1-nu}\\
                 =        &\, \frac{nu+u(1-nu)+nu^2}{1-nu} = \frac{(n+1)u}{1-nu} \leqslant \frac{(n+1)u}{1-(n+1)u} = \gamma_{n+1}.
\end{align*}
The case $\rho_{n+1}=-1$ can be proved, similarly.
$\qed$
\vsp

\subsection{Multiplication}
Assume that \( x_1, x_2, \ldots, x_n \in \mathbb{F} \) are nonzero numbers and let \( s_n = x_1 x_2 \cdots x_n \) be their product. We assume that the multiplication is carried out from left to right.The following simple Python function demonstrates how this operation is performed.
\begin{shaded}
\vspace*{-0.3cm}
\begin{verbatim}
def Multiplication(x):
    s = x[0]
    for k in range (1,len(x)):
        s  = s * x[k]
    return s
\end{verbatim}
\vspace*{-0.3cm}
\end{shaded}
\noindent
In this function, \verb+x+ is a list of numbers, and the product is computed iteratively by multiplying each element of the list from left to right. The variable s holds the intermediate product at each step, and the final product is returned at the end.

To analyse this algorithm, let \( s_k = x_1 x_2 \cdots x_k \) denote the \( k \)-th partial product.
In the following analysis, and throughout the lecture, a letter with hat (e.g. \(\wh x\) ) denotes a computed (floating-point) quantity.
 Using the standard model
\eqref{fl_op_rule} we have
\begin{align*}
\sh_1 &= \fl(x_1) = x_1 \quad (\mbox{indeed }\delta_1 = 0)\\
\sh_2 &= \fl(\sh_1x_2) = \sh_1x_2(1+\delta_2) = x_1x_2(1+\delta_2)\\
\sh_3 &= \fl(\sh_2x_3) = \sh_2x_3(1+\delta_3) = x_1x_2(1+\delta_2)x_3(1+\delta_3).
\end{align*}
for $|\delta_k|\leqslant u$, $k=2,3$. Continuing in the same way, we obtain
\begin{align}\label{prod_backward}
\sh_n  = x_1x_2(1+\delta_2)x_3(1+\delta_3)\cdots x_n (1+\delta_{n}), \quad |\delta_k|\leqslant u,
\end{align}
and using Lemma \ref{lem:gamma_n} we can write
\begin{equation}\label{prod_forward}
\begin{split}
\sh_n  &= x_1x_2\cdots x_n (1+\delta_2)(1+\delta_3)\cdots(1+\delta_{n})\\
& = s_n (1+\theta_{n-1}), \quad  |\theta_{n-1}|\leqslant \gamma_{n-1}.
\end{split}
\end{equation}
The expression \eqref{prod_backward} is a {\em backward error} analysis and can be interpreted as follows: the computed
product 
$$
\sh_n = \wh x_1 \wh x_2 \cdots \wh x_n
$$ 
is indeed the exact product of perturbed data
$$
\xh_1=x_1, \quad \xh_2=x_2(1+\delta_2),\; \ldots, \; \xh_n=x_n(1+\delta_n), \quad |\delta_k|\leqslant u.
$$
Each relative perturbation is bounded by $u$,
so the perturbations are tiny. This shows that the computed solution $\sh_n$ is the product of {\em nearby} data $\wh x_k$.  

On the other hand, \eqref{prod_forward} gives the {\em forward error} bound
\vsp 
$$
\frac{|\sh_n-s_n|}{|s_n|}\leqslant \gamma_{n-1},
$$
which estimates the difference between the computed solution $\sh_n$ and the exact solution $s_n$.

The above analysis exhibits that multiplication algorithm is {\em backward stable} independent of $n$ (the number of operands), and {\em forward stable} if $n$ (number of operands) is not a large number.
Note than $\gamma_n$ is increasing in $n$. 
\vsp
\subsection{Summation}
 Assume \( x_1, \ldots, x_n \in \mathbb{F} \) and let \( s_n = x_1 + \cdots + x_n \). If the sum is computed using the usual recursive summation, as the code below,
 \begin{shaded}
\vspace*{-0.3cm}
\begin{verbatim}
def Summation(x):
    s = x[0]
    for k in range (1,len(x)):
        s  = s + x[k]
    return s
\end{verbatim}
\vspace*{-0.3cm}
\end{shaded}
\noindent
 then after some detailed analysis similar to what we did in the previous example for multiplication, we obtain
\begin{equation}\label{sum_backward}
  \sh_n = x_1(1+\theta'_{n-1}) + x_2(1+\theta_{n-1})+x_3(1+\theta_{n-2})+\cdots + x_n(1+\theta_1), \quad |\theta_k|\leqslant \gamma_k,
\end{equation}
and
\begin{equation}\label{sum_forward}
  |\sh_n -s_n| \leqslant  |x_1| \gamma_{n-1} + |x_2|\gamma_{n-1}+|x_3|\gamma_{n-2}+\cdots + |x_n|\gamma_1.
\end{equation}
From \eqref{sum_backward} we observe that the computed solution $\sh_n$ is indeed the exact sum of nearby data
$$
\xh_1 = x_1(1+\theta'_{n-1}), \quad \xh_k = x_k(1+\theta_{n-k+1}),\mbox{  with  }|\theta'_{n-1}|\leqslant \gamma_{n-1},\; |\theta_{n-k+1}|\leqslant \gamma_{n-k+1}.
$$
This is a backward error analysis. Although at first glance summation might appear to result in smaller roundoff errors than multiplication, a comparison with the multiplication algorithm reveals that the recursive summation algorithm has larger backward error bounds. This is because the backward errors \(u\) for multiplication are replaced by values \(\gamma_{n-1}\) and \(\gamma_{n-k+1}\) in the summation algorithm.
From \eqref{sum_forward}, we have the forward error bound
\begin{equation}\label{sum_forward2}
  |\sh_n - s_n| \leqslant \gamma_{n-1} \sum_{k=1}^{n} |x_k|.
\end{equation}
This upper bound holds independently of the summation order. However, the upper bound \eqref{sum_forward} can be minimized if the terms \(x_k\) are added to the sum in increasing order of magnitude. By summing smaller terms first, the impact of roundoff errors is reduced, leading to a more accurate result. This strategy leverages the associative property of addition to mitigate the accumulation of roundoff errors.

In the above analysis, we observed that summation in floating-point arithmetic can introduce significant errors in the final computed solution. To address this issue, various techniques have been developed for more accurate summation. One such clever algorithm is {\em compensated summation}, a method that incorporates a correction term to reduce rounding errors. This technique captures the rounding errors and feeds them back into the summation process, thereby improving accuracy \cite{Kahan:1965}.
For more details on compensated summation and other advanced techniques, see \cite{Higham:2002}.
\vsp 

\subsection{Inner product and matrix multiplications}

Consider the inner product $s_n = x^Ty$, where $x,y\in\F^n$.
We assume that the evaluation of $s_n=x_1y_1+\cdots+x_ny_n$ is performed from left to right.
Let $s_k=x_1y_1+\cdots+x_ky_k$ denote the $k$-th partial sum. Using the standard floating-point model \eqref{fl_op_rule} and without fused multiply-add (FMA) operations, we have, after some calculations
\begin{equation*}
  \sh_n = x_1y_1(1+\theta'_n) + x_2y_2(1+\theta_n) + x_3y_3(1+\theta_{n-1})+\cdots + x_ny_n(1+\theta_2)
\end{equation*}
for $|\theta'_{n}|\leqslant \gamma_n$ and $|\theta_{k}|\leqslant \gamma_k$ for $k=2,\ldots,n-1$. This is a backward error analysis, showing that
$$
\sh_n = \wh x_1\wh y_1+\cdots+\wh x_n\wh y_n
$$
where
$$
\xh_k = x_k,\; k=1,\ldots,n, \quad \yh_1=y_1(1+\theta'_n), \quad \yh_k = y_k(1+\theta_{n-k+2}),\; k=2,\ldots,n.
$$
Alternatively, we could perturb $x_k$ and leave $y_k$ alone. Using a vector form, we can write
\begin{equation}\label{inner_backward}
\sh_n = \fl(x^Ty) = x^T(y+\delta y) = (x+\delta x)^Ty,\quad |\delta x|\leqslant \gamma_n|x|,\quad |\delta y|\leqslant \gamma_n|y|
\end{equation}
with $|x|$ denoting the vector with elements $|x_k|$. Inequalities between vectors
(and, later, matrices) are understood componentwise.  A forward error bound follows from \eqref{inner_backward}:
\begin{equation}\label{inner_backward2}
|x^Ty -\fl(x^Ty)| \leqslant \gamma_n \sum_{k=1}^{n}|x_ky_k|=\gamma_n |x|^T|y|.
\end{equation}
If $y = x$, a high
relative accuracy is obtained for computing $x^Tx$. However, in general, high relative accuracy is not
guaranteed if $|x^Ty|\ll |x|^T|y|$.
\vsp 

\begin{workout}
To compute the inner product $s_n=x^Ty$ of vectors $x,y\in\F^n$ for an even $n$, assume $m=n/2$, $s_1 = x(1:m)^Ty(1:m)$, $s_2=x(m+1:n)^Ty(m+1:n)$, and $s_n= s_1+s_2$. Prove the forward error bound
$$
|\sh_n-s_n|\leqslant \gamma_{n/2+1}|x|^T|y|
$$
which shows
the error bound is almost halved by separating the inner product in two pieces. Generalized the idea by
breaking the inner product into $k$ pieces and find the optimal $k$.
\end{workout}
\vsp 

With the analysis of the inner product in hand, it becomes straightforward to analyze matrix-vector and matrix-matrix multiplications. Consider
$A\in \F^{m\times n}$, $x\in\F^n$ and $y=Ax$. The vector $y$ can be formed by computing $m$ inner products $y_k = a_{\cdot k}^Tx$, $k=1,\ldots,m$, where
$a_{\cdot k}^T$ is the $k$-th row of $A$.  From \eqref{inner_backward2} we have
$$
\wh y_k = (a_{\cdot k}+\delta a_{\cdot k})^Tx, \quad |\delta a_{\cdot k}|\leqslant \gamma_n|a_{\cdot k}|, \quad k=1,\ldots, m,
$$
which gives the backward error
\begin{equation}\label{mat-vec:backward}
\yh = (A+\delta A)x, \quad |\delta A|\leqslant \gamma_n |A|.
\end{equation}
From \eqref{mat-vec:backward},  we obtain a forward error bound as
\begin{equation}\label{mat-vex:forward}
|\yh-y|\leqslant \gamma_n |A||x|.
\end{equation}
Note that comparison operators (equalities and inequalities) between matrices and vectors are understood component-wise. Normwise bounds readily follow. For example, it is not difficult to show that if $|x|\leqslant |y|$, then $\|x\|_p\leqslant \|y\|_p$ for $p=1,2,\infty$, and indeed for any $p$. The same holds true for matrices: if $|A|\leqslant |B|$, then $\|A\|_p\leqslant \|B\|_p$ but just for $p=1,\infty$. For $p=2$ we have $\|A\|_2\leqslant \sqrt{\min\{m,n\}}\|B\|_2$. Thus,
the component-wise forward bound \eqref{mat-vex:forward} results in
$$
\|\yh-y\|_p\leqslant \gamma_n \|A\|_p \|x\|_p, \quad p=1,\infty.
$$
and
$$
\|\yh-y\|_2\leqslant \sqrt{\min\{m,n\}} \gamma_n \|A\|_2 \|x\|_2.
$$
\vsp
\begin{workout}
Assume that $A\in \F^{m\times p}$ and $B\in\F^{p\times n}$ and $C=AB$. Derive forward and backward component-wise error bounds for this matrix-matrix multiplication. 
\end{workout}
\vsp

\subsection{Cancellation}

In this section, we discuss cancellation, a particularly dangerous phenomenon in numerical computing that can cause significant accuracy loss during simple operations like subtracting machine numbers with the same sign that have been previously rounded and share several leading digits. The same issue occurs in floating-point addition of machine numbers with opposite signs that have similar magnitudes and have been previously rounded.

To be precise, consider a (potentially long) chain of numerical computations where \(\wh x\) and \(\wh y\) are two machine numbers that are our existing \(p\)-digit approximations to some real quantities \(x\) and \(y\), respectively. These quantities \(x\) and \(y\) may have infinitely many digits. Assume that among the significant digits of \(\wh x\) and \(\wh y\), \(k\) of the leading digits are the same. Now, consider the task of computing \(\fl(\wh x - \wh y)\). Ideally, we would want to compute the exact value \(x - y\), but in practice, we can only compute \(\fl(\wh x - \wh y)\) with \(p\)-digit precision. We have
\vsp 
\begin{align*}
 & (d_0.d_1 d_2 \cdots d_{k-1}\, d_k\, d_{k+1} \cdots d_{p-1})_{\beta} \times \beta^{e} \\
-\ &  (d_0.d_1 d_2 \cdots d_{k-1}\; c_k\; c_{k+1} \cdots c_{p-1})_{\beta}\times \beta^{e}\\
 \cline{1-2}
=\ &\  (0\ .\ 0\  0\ \cdots\ 0\ \ f_k\ f_{k+1} \cdots f_{p-1})_{\beta}\times \beta^{e}
\end{align*}
When normalized, this result is
\[
(f_k. f_{k+1} \cdots f_{p-1}\underbrace{00 \cdots 0}_{k \mbox{ \footnotesize{times}}})_{\beta}\times \beta^{e-k}.
\]
Here, \(k\) last significant digits in the mantissa become zero, and their truly correct values are not known. This loss of significant digits in the subtraction of two close machine numbers is known as {\bf cancellation error}.

To understand why cancellation in the operation \(\fl(\wh x - \wh y)\) is dangerous, consider that the resulting digits \((f_k. f_{k+1} \cdots f_{p-1})_{\beta}\) have been computed by subtracting the \((d_k\, d_{k+1} \cdots d_{p-1})_{\beta}\) part of \(\wh x\) and the \((c_k\, c_{k+1} \cdots c_{p-1})_{\beta}\) part of \(\wh y\). These tail digits of \(\wh x\) and \(\wh y\) are the most affected by rounding errors, and thus, have a higher chance of being different from the corresponding digits of \(x\) and \(y\).
As a result, cancellation can cause a significant loss of accuracy because the number of correct digits in the output \(\fl(\wh x - \wh y)\) might be much less than the number of correct digits in the inputs \(\wh x\) and \(\wh y\). This is why the cancellation phenomenon is sometimes referred to as {\em catastrophic cancellation}.
\vsp

\begin{example}
We assume that
\begin{align*}
& x = 1.23456702645\\
& y = 1.2345664932563685
\end{align*}
and compute $x - y$ in a floating point system with $\beta=10$ and $p=7$. 
Notice that in a decimal system with $p=7$ we basically have seven significant decimal digits which is similar to using single precision in the IEEE standard. With rounding to nearest we have
\begin{align*}
& \wh x = \fl(x) = 1.234567\\
& \wh y = \fl(y) = 1.234566.
\end{align*}
Therefore, we can write 
\[
\fl(\wh x - \wh y) = \fl(\fl(x) - \fl(y)) =\fl(0.000001) = 1.000000 \times 10^{-6}.
\]
On the other hand, the correct result $x - y$, which in general can be computed only in theory, is as follows:
\begin{align*}
x - y &= 1.2345670264500000 - 1.2345664932563685\\
& = 0.0000005331936315 = 5.331936315 \times 10^{-6}.
\end{align*}
The relative error then can be obtained as 
\vsp 
\[
\frac{|\fl(\wh x - \wh y) - (x-y)|}{|x-y|} = 
\frac{|1.000000 \times 10^{-6} - 5.331936315 \times 10^{-6}|}{|5.331936315 \times 10^{-6}|} \approx 0.87
\]
which is huge in comparison with the corresponding unit roundoff ($ \approx 10^{-7}$). Our inputs $\wh x$ and $\wh y$ had seven correct significant digits while our output $\fl(\wh x - \wh y)$ does not have even one significant correct digit!
\end{example}
\vsp
\begin{example}
Let us design an algorithm to approximate the irrational number $\pi$ by considering it as the circumference of a semi-circle with radius $1$. To this aim, we divide the semi-circle to $n$ equi-length arcs to obtain an internal regular semi-polygons;
see Figure \ref{fig:polygons}.

\begin{center}
\includegraphics[scale=0.6]{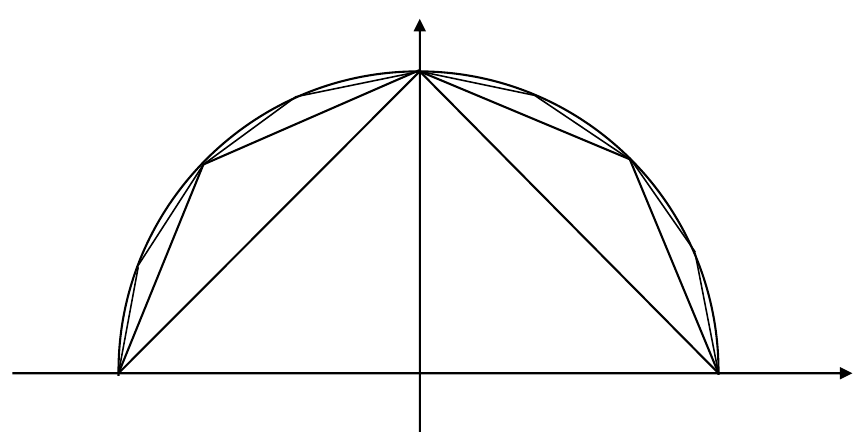}
\captionof{figure}{Polygons approximation of a circle.}
\label{fig:polygons}
\end{center}

The circumference of the semi-circle can be approximated by the sum of lengths of its line segments (circumference of the semi-polygons). Since the length of each side of the polygons is $\sin(\pi/n)$, the circumference of the semi-polygons is equal to 
$$
f(n) = n\sin\frac{\pi}{n}.
$$
Obviously, $\displaystyle\lim_{n\to\infty}f(n) = \pi$. Let us increase $n$ by powers of $2$, i.e., $n=1,2,4,\ldots,2^k,\ldots$, and define
\begin{equation*}
  p_{k-1} = f(2^{k}) = 2^k\sin 2^{-k}\pi, \quad k=1,2,\ldots
\end{equation*}
From identity 
$$
\sin\frac{\alpha}{2} = \sqrt{\tfrac{1}{2}(1-\cos \alpha)},\quad 0\leqslant\alpha\leqslant 2\pi,
$$
we can write for $\alpha=2^{-k}\pi$,
$$
p_{k}=2^{k+1}\sin \frac{\alpha}{2} = 2^{k+1}\sqrt{\tfrac{1}{2}(1-\cos \alpha)}=2^{k+1}\sqrt{\frac{1}{2}\left(1-\sqrt{1-\sin^2 \alpha}\right)}.
$$
Using the fact that $\sin \alpha = 2^{-k}p_{k-1}$, we have the following recursive formula for the $\{p_k\}$ sequence:
\begin{equation}\label{pk-seq}
  p_{k}=2^{k+1}\sqrt{\frac{1}{2}\left(1-\sqrt{1-[2^{-k}p_{k-1}]^2}\right)}, \quad k=1,2,\ldots, \quad p_0=2.
  \vsp
\end{equation}
Theoretically, the sequence $\{p_k\}$
tends to $\pi$ as $k$ increases. Let's see what happens if we compute it numerically with a Python code.
\begin{shaded}
\vspace*{-0.3cm}
\begin{verbatim}
import numpy as np
import matplotlib.pyplot as plt
K = 30
p = np.zeros(K)
p[0] = 2
for k in range (1,K):
    p[k] = 2**(k+1)*np.sqrt(0.5*(1-np.sqrt(1-2**(-2*k)*p[k-1]**2)))
plt.figure()
plt.semilogy(np.arange(K), abs(np.pi-p), marker = 's', color ='red')
plt.title("Approximation of $\pi$")
plt.xlabel("$k$"), plt.ylabel("$|\pi-p_k|$"),
\end{verbatim}
\vspace*{-0.3cm}
\end{shaded}
\noindent

The error \(|\pi - p_k|\) is shown on the left-hand side of Figure \ref{fig:pi_app}. Initially, the error decreases to about \(10^{-9}\) until \(k = 14\), but unexpectedly begins to increase for larger \(k\) values, up to \(k=29\). We observe a $V$-shape error plot that reveals an instability in numerical computation.

\begin{center}
\includegraphics[scale=0.6]{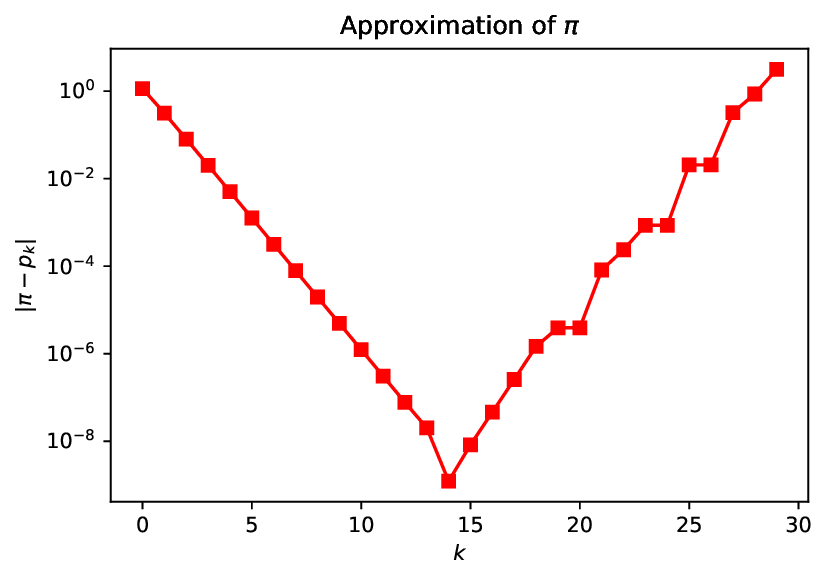}\includegraphics[scale=0.6]{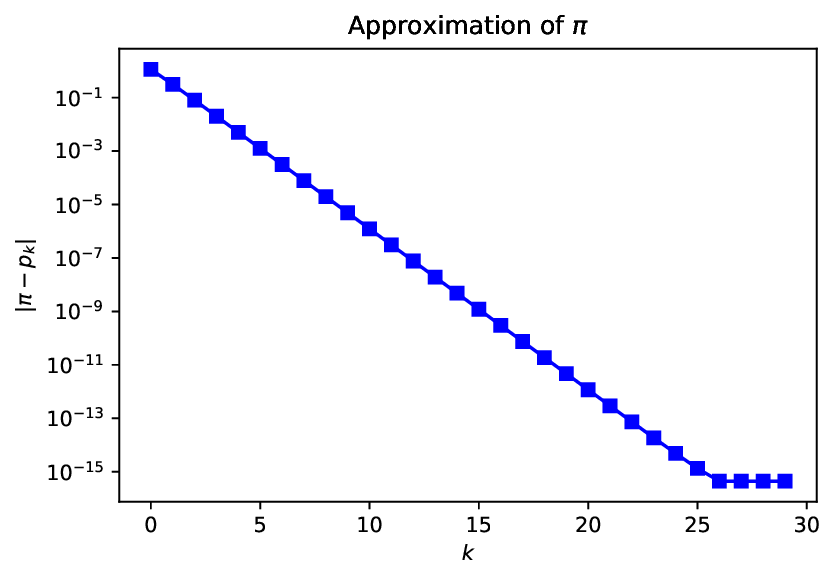}
\captionof{figure}{Approximation of $\pi$ with two mathematically equivalent, but numerically different, formulas \eqref{pk-seq} (left) and \eqref{pk-seq2} (right).}
\label{fig:pi_app}
\end{center}

Let us reformulate \eqref{pk-seq} using the identity
\begin{equation}\label{identity_1x2}
1-x=\frac{1-x^2}{1+x}, \quad x\neq -1
\end{equation}
for $x =\sqrt{1-[2^{-k}p_{k-1}]^2}$.  The new recursive sequence, which is theoretically equivalent to \eqref{pk-seq}, is
\begin{equation}\label{pk-seq2}
  p_{k}=p_{k-1}\sqrt{\frac{2}{1+\sqrt{1-[2^{-k}p_{k-1}]^2}}}, \quad k=1,2,\ldots,\quad  p_0=2.
  \vsp
\end{equation}
We update our Python code by replacing the line inside the loop with the new formula \eqref{pk-seq2}. The resulting error is plotted on the right-hand side of Figure \ref{fig:pi_app}. Notably, the error function now monotonically decreases to the level of machine precision. For a clearer comparison, both graphs are presented together in Figure \ref{fig:pi_app2}.

\begin{center}
\includegraphics[scale=0.75]{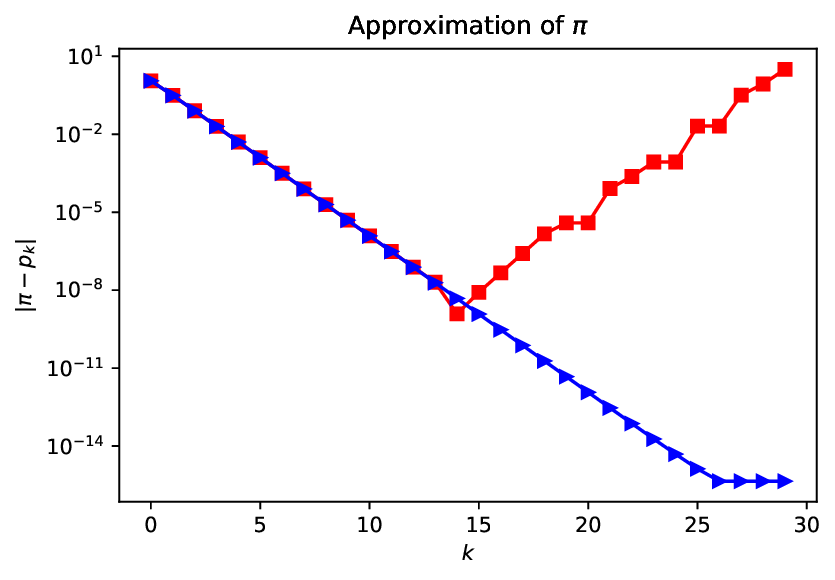}
\captionof{figure}{Approximation of $\pi$ with two mathematically equivalent, but numerically different, formulas \eqref{pk-seq} (red line and square markers) and \eqref{pk-seq2} (blue line and triangle markers).}
\label{fig:pi_app2}
\end{center}

Formulas \eqref{pk-seq} and \eqref{pk-seq2} are mathematically equivalent, but there exists a source of cancellation error in formula \eqref{pk-seq} that makes it improper for floating-point arithmetic. Consider the term $1 - \sqrt{1 - [2^{-k}p_{k-1}]^2}$. As $k$ increases, $2^{-k}p_{k-1}$ tends to zero, causing $\sqrt{1 - [2^{-k}p_{k-1}]^2}$ to approach $1$. Consequently, there is a risk of  cancellation error in the subtraction $1 - \sqrt{1 - [2^{-k}p_{k-1}]^2}$ for large values of $k$. In formula \eqref{pk-seq2}, we deliberately avoid this subtraction by using the simple identity \eqref{identity_1x2}. This minor adjustment eliminates the cancellation error and ensures the stability of the algorithm.
\end{example}
\vsp
In numerical calculations, if possible one should try to avoid formulas with subtraction of close floating point numbers that give rise
to cancellation error.
\vsp
\begin{remark}
Beyond an illustration for cancelation error, the above examples reveals a
fundamental insight in numerical computations:
``mathematically equivalent'' formulas or algorithms are not in general ``numerically
equivalent''.
\end{remark}
\vsp 
\begin{workout}
Reformulate the following expressions to avoid possible cancellations in computation.
\begin{align*}
  &1-\cos x, \quad |x|\ll 1 \\
  &\sin x-\cos x ,\quad |x|\approx \frac{\pi}{4} \\
  &\ln \left(\sqrt{1+x^2}-x  \right) ,\quad |x|\gg 1
\end{align*}
\end{workout}
\vsp
We can apply our knowledge of the condition number to observe that subtraction of two close numbers (cancellation) is dangerous. This can be done in various ways but for our purposes it is enough to simply allow just one of the inputs $x$ and $y$ to vary and assume that the other one is constant. In particular, we consider the condition number of the problem of evaluating the function
$f(x) = x - c$ for variable $x$ and constant $c$. According to definition of condition number of univariate functions, i.e. \eqref{cond:def1}, we have
\vsp
\[
(\cond f)(x) = \frac{|x| \cdot 1}{|x-c|} =  \frac{|x|}{|x-c|} .
\]
This means that when $x$ and $c$ are close to each other, the denominator of the condition number becomes small, leading to a large condition number. Therefore, it is not surprising that the cancellation phenomenon is dangerous. It is analogous to solving an ill-conditioned problem, where small changes in the input can result in significant changes in the output.
\vsp 

\section{Complexity of an algorithm}
The {\bf complexity} (or {\bf computational cost}) of an algorithm is its executing time and the amount of resources required to run it.
Particular focus is given to {\em time} and {\em space} (memory) requirements.
Calculating the complexity of an algorithm is therefore a part of the
analysis of its efficiency.

\begin{definition}
{\em Time complexity} of an algorithm is generally expressed as the total number of required {\bf fl}oating-point {\bf op}eration{\bf s} ({\bf flops}) $\{+,-,\times,/\}$ on its input data to produce the final output.
\end{definition}

If the input data is of size $n$, then the time complexity will be a function of $n$, say $f(n)$.
Usually, we focus on the behavior of the complexity for large $n$, that is on its asymptotic behavior when $n$ tends to the infinity. Therefore, the complexity is generally expressed by using big $\mathcal O$ notation.
\vsp

\begin{example}
The inner product of two vectors
$x,y\in\R^n$ that is $s = x_1y_1+x_2y_2+\cdots + x_ny_n$ requires $n$ multiplications and $n-1$ additions. The complexity for this computation is $f(n)=2n-1$ flops which is asymptotically expressed by $\mathcal O(n)$ flops.
\end{example}
\vsp

\begin{example}
The time complexity for matrix-vector product $Ax$ for $A\in\R^{n\times n}$ and $x\in \R^n$ is $f(n)=2n^2-n$ flops (why?) or asymptotically $\mathcal O(n^2)$ flops. The time complexity for matrix-matrix product $AB$ for $A,B\in \R^{n\times n}$ is $f(n)=2n^3-n^2$ flops (why?) or $\mathcal O(n^3)$ flops.
\end{example}
\vsp 

Sometimes, the leading coefficient in $f(n)$ is also mentioned in the asymptotic expression of the time complexity. For example, the time complexity of the matrix-matrix product is reported as $2n^3$ flops. The reason should be clear; there is a huge difference between
$1000n^3$ and $0.1n^3$, even for large values of $n$.
\vsp

\begin{example}
The complexity of Gauss elimination algorithm for solving the linear system $Ax=b$ for $n\times n$ matrix $A$ is
$\frac{2}{3}n^3$ flops.
\end{example}
\vsp
The usual units of time (seconds, minutes etc.) are not used for time complexity because they are dependent on the choice of a specific computer and on the evolution of technology.
\vsp

\begin{definition}
{\em Space complexity} is generally expressed as the amount of memory required by an algorithm to load the inputs, execute, and produce the final solution (output).
\end{definition}
\vsp

The amount of memory needed depends on a variety of things such as the programming language, the compiler, or even the machine running the algorithm. Here, we just follow a simple rule:
if we need to create a scaler this will require $1$ {\bf u}nit {\bf o}f {\bf s}pace ({\bf uos}); if we create an array of size $n$, this will require $n$ uos; and if we create a matrix (a two-dimensional array) of size $m\times n$, it will require $mn$ uos.
\
\begin{example}
The space complexity of the standard algorithm for computing the inner product of two vectors $x,y\in\R^n$ is $2n$ uos to store
$x$ and $y$, and $n$ auxiliary uos to store partial summations. If we let partial summations overwrite, only $1$ auxiliary uos is needed. The output can also be stored in that auxiliary space. Consequently, the space complexity for this algorithm is $2n+1$ uos.
\end{example}
\
\begin{example}
Consider the Gauss elimination algorithm for solving the system $Ax=b$ for square matrix $A$ of size $n\times n$ and vector $b\in \R^n$. The input space is clearly $n^2+n$ uos. The auxiliary space needed depends on how much thrifty the algorithm is implemented. Remember that the process of Gauss elimination consists of $n-1$ steps to convert $A$ to an upper triangular matrix.
See the following illustration for $n=4$.

\begin{equation*}
\begin{array}{ccccccc}
  \begin{bmatrix}
    \times &\times &\times & \times \\
    \times &\times &\times & \times \\
    \times &\times &\times & \times \\
    \times &\times &\times & \times
  \end{bmatrix}&\longrightarrow &
  \begin{bmatrix}
    \times &\times &\times & \times \\
    0 &\times &\times & \times \\
    0 &\times &\times & \times \\
    0 &\times &\times & \times
  \end{bmatrix}&\longrightarrow &
  \begin{bmatrix}
    \times &\times &\times & \times \\
    0 &\times &\times & \times \\
    0 &0 &\times & \times \\
    0 &0 &\times & \times
  \end{bmatrix}&\longrightarrow &
  \begin{bmatrix}
    \times &\times &\times & \times \\
    0 &\times &\times & \times \\
    0 &0 &\times & \times \\
    0 &0 &0 & \times
  \end{bmatrix}\\
  A&&A^{(1)}&&A^{(2)}&& A^{(3)}
  \end{array}
\end{equation*}

In each step $k$,
the same operations are applied on vector $b^{(k-1)}$ to convert it to the new vector $b^{(k)}$. Finally, the equivalent upper triangular system
$$
A^{(n-1)}x = b^{(n-1)}
$$
should be solved using a backward substitution for solution $x$.

The process of elimination can be carried out without any auxiliary space if we overwrite $A$ and $b$ in each step, and store the auxiliary variables in the course of elimination in the lower diagonal part of $A$ (zero positions). The backward substitution
needs $n$ auxiliary uos to store the output $x$. Consequently, the total space complexity is $n^2+2n$ uos.

This complexity analysis is valid for Gauss elimination algorithm without pivoting. For Gauss elimination with pivoting, a $2n+1$ extra uos is required, $n+1$ uos to interchange the rows of $A$ and $b$ in each step, and $n$ uos to keep the permutation history in an array of size $n$.
\end{example}
\vsp

Generally, when ``complexity'' or ``computational cost'' is used without being further specified, we mean the worst-case time complexity of the algorithm.
\vsp

\section*{Additional workouts}
Here are some additional exercises to help solidify your understanding of the concepts covered in the lecture. Some of these exercises are adapted from the references mentioned at the beginning of the lecture.
\vsp
\begin{workout}
In the IEEE standard with double precision determine an interval for which the distance between the floating-point numbers is exactly $1$. 
\end{workout}
\begin{workout}

How many double precision floating-point numbers do exist between two consecutive nonzero single precision floating-point numbers? 
\end{workout}

\begin{workout}
Show that
$$
0.1 = \sum_{k=1}^{\infty} \left(2^{-4k}+2^{-4k-1}\right)
$$
and conclude $(0.1)_{10} =(0.000\overline{1100})_2 $. The last four binary digits are repeated. In the IEEE standard with single precision show that if $\wh x = \fl(0.1)$ then
$$
\frac{x-\wh x}{x}=\frac{1}{4}u.
$$
\end{workout}
\begin{workout}
Let $x$ be a floating point number in IEEE double
precision arithmetic satisfying $1 \leqslant x < 2$. Show that $\fl(x \times (1/x))$ is either $1$ or
$1 -\ep_M/2$.
\end{workout}

\begin{workout}
Show by an example that the inequalities $x \leqslant \fl((x + y)/2) \leqslant y$, where $x$ and $y$ are floating
point numbers with $x\leqslant y$, can be violated in base $10$ arithmetic. Show that
$x \leqslant \fl(x + (y - x)/2) \leqslant y$ hold true in any base $\beta$ arithmetic.
\end{workout}

\begin{workout}
Show that with gradual underflow, if $x$ and $y$ are
floating-point numbers and $\fl(x \pm y)$ underflows then $\fl(x\pm y) = x \pm y$. (no roundoff error for addition and subtraction when result falls in the underflow range).
\end{workout}

\begin{workout}
(Kahan) Let $A = \begin{bmatrix}
                    a & b \\
                    c & d
                  \end{bmatrix}$.
Show that with the use of a fused multiply-add operation the
algorithm
\begin{align*}
w &= bc\\
e &= w - bc\\
x &= (ad - w) + e
\end{align*}
computes $x = \mathrm{det}(A)$ with high relative accuracy.
\end{workout} 
 
\begin{workout}
A part of this course is devoted to numerical solution of ordinary differential equations (ODEs). 
The ODE is posed on a specific interval, $[a,b]$, and a numerical method works with a discretization of this interval with equidistance points
$x_k:=a+hk$, $k=0,1,\ldots,n$ with $h=(b-a)/n$. Compare the accuracy of the following formulas for computing $x_k$. Assume that $a$ and $b$ are floating-point numbers but $h$ is not. 
\begin{align*}
  \mathrm{(I)}\quad & x_k = x_{k-1}+h \\
  \mathrm{(II)}\quad & x_k=a+kh  \\
  \mathrm{(III)}\quad  & x_k = a(1-k/n)+(k/n)b
\end{align*}
Usually, without any reason, (I) is used by users. 
\end{workout}
 
  \begin{workout}
Compute the condition number of the following functions, and discuss any
possible ill-conditioning.
\begin{align*}
  & f(x) = \sin^{-1}x \\
  & f(x) = x^{1/n},\quad x>0,\quad n\in\N \\
  & f(x) = \cos (x),\quad |x|\leqslant \frac{\pi}{2}
\end{align*}
\end{workout}
 
\begin{workout}
Show that 
$$
(\cond\, fg)(x)\leqslant (\cond\, f)(x)+(\cond\, g)(x).
$$
What can be said
about $(\cond\, f/g)(x)$?
\end{workout}
  
\begin{workout}
Assume that $h(t)=g(f(t))$. Compute the condition number of
$h$ in terms of the condition of $g$ and $f$.
Apply it on function $h(t) = \frac{1+\sin t}{1-\sin t}$ to compute its condition number at $t=\pi/4$. 
\end{workout}
 
\begin{workout}
Consider the algebraic equation
$$
x^n +ax-1=0,\quad a>0,\quad n\geqslant 2.
$$
Show that the equation has exactly one positive root $\xi(a)$.
Show that $\xi$ is well-conditioned with respect to perturbations in $a$. 
\end{workout}
 
\begin{workout}
For a nonsingular matrix $A\in \R^{n\times n}$ show that 
$$
\cond_2(A)=\sqrt{\frac{\lambda_{\max}(A^TA)}{\lambda_{\min}(A^TA)}}.
$$
If $A$ is further symmetric then show 
$$
\cond_2(A)=\frac{\max|\lambda_k(A)|}{\min|\lambda_k(A)|}.
$$
\end{workout}
 
\begin{workout}
Assume that $A\in \R^{n\times n}$ and $\lambda$ is an eigenvalue of $A^TA$. Show that
$$
0\leqslant\lambda \leqslant \|A^T\| \cdot \|A\|
$$
where $\|\cdot\|$ is an operator norm. Using this show that 
$$
[\cond_2(A)]^2\leqslant \cond_1(A)\cdot \cond_\infty(A).
$$
provided that $A$ is nonsingular.
\end{workout}
\vsp


\begin{thebibliography}{9}


\bibitem{Dahlquist-Bjork:2008}
G. Dahlquist, \AA. Bj\"{o}rck,  {\em Numerical Methods in Scientific Computing}, Volume 1, SIAM, Philadelphia, PA, 2008.

\bibitem{Demmel:1984}
J. W. Demmel, Underflow and the reliability of numerical software. SIAM J. Sci. Stat.
Comput., 5(4):887--919, 1984.

\bibitem{Heath:2018} M. T. Heath, {\em Scientific Computing, an Introductory Survey}, revised 2nd edition, SIAM, Philadelphia, PA, 2018.

\bibitem{Gautschi:2012} W. Gautschi, {\em Numerical Analysis}, 2nd edition, Springer, 2012.

\bibitem{Higham:2002} Nicholas J. Higham, {\em Accuracy and Stability of Numerical Algorithms}, 2nd edition, SIAM, Philadelphia, PA,
2002.

\bibitem{Kahan:1965} W. Kahan, Further remarks on reducing truncation errors, Comm. ACM, 8(1) (1965) 40.


\end{thebibliography}
\end{document}